# STABLE ERGODICITY OF CERTAIN LINEAR AUTOMORPHISMS OF THE TORUS

FEDERICO RODRIGUEZ HERTZ

ABSTRACT. We prove that some ergodic linear automorphisms of $\mathbb{T}^N$ are stably ergodic, i.e. any small perturbation remains ergodic. The class of linear automorphisms we deal with includes all non-Anosov ergodic automorphisms when $N = 4$ and so, as a corollary, we get that every ergodic linear automorphism of $\mathbb{T}^N$ is stably ergodic when $N \leq 5$.

## 1. Introduction

Given a matrix $A \in SL(N, \mathbb{Z})$ we define the linear automorphism of $\mathbb{T}^N = \mathbb{R}^N / \mathbb{Z}^N$ generated by $A$, which we also denote $A$, via the following diagram

$$
\begin{array}{ccc}
\mathbb{R}^N & \xrightarrow{A} & \mathbb{R}^N \\
\downarrow & & \downarrow \\
\mathbb{T}^N & \xrightarrow{A} & \mathbb{T}^N
\end{array}
$$

where the down arrows are the canonical projections.

That a linear automorphism of the torus is ergodic if and only if it has no eigenvalue being a root of the unity is known since the work of Halmos [Ha]. Only after the work of Anosov [A] it turns out that some linear automorphisms (the so called Anosov linear automorphisms) were in fact stably ergodic. In this case, stable ergodicity comes from the stability of Anosov diffeomorphisms and its ergodicity. Around 1969, Pugh and Shub began to study the stable ergodicity of diffeomorphisms and they wondered if the following linear automorphism of $\mathbb{T}^4$

$$
\begin{pmatrix}
0 & 0 & 0 & -1 \\
1 & 0 & 0 & 8 \\
0 & 1 & 0 & -6 \\
0 & 0 & 1 & 8
\end{pmatrix}
$$

was stably ergodic. Moreover, in [HPS], Hirsh, Pugh and Shub posed the following question:

**Question 1.** *Is any ergodic linear automorphism of $\mathbb{T}^N$ stably ergodic?*

The goal of this paper is to give a positive answer to this question for some linear automorphisms. We say that a linear automorphism of $\mathbb{T}^N$ is pseudo-Anosov if it has no eigenvalue being root of the unity, its characteristic polynomial irreducible over the integers and not a polynomial in $t^n$ for some $n \geq 2$. Let us explain why the name pseudo-Anosov. Take a surface $S$ of genus $g$, a homeomorphism $h : S \to S$ and look at the action $h_* : H_1(S, \mathbb{Z}) \to H_1(S, \mathbb{Z})$ of $h$ on the first homology group of the surface. As $H_1(S, \mathbb{Z}) \simeq \mathbb{Z}^{2g}$, $h_*$ induces a matrix in $SL(2g, \mathbb{Z})$ and hence a linear automorphism of $\mathbb{T}^{2g}$. We have that

*Date*: October 30, 2018.
This work has been partially supported by IMPA/CNPq.

if $h_*$ induces a pseudo-Anosov linear automorphism on $\mathbb{T}^{2g}$ then $h$ is isotopic to a pseudo-Anosov homeomorphism of $S$. In order to state our results let us denote the center space by $E^c$. (i.e the eigenspace of the eigenvalues with modulus one)

**Theorem 1.1.** *For $N \geq 6$, any pseudo-Anosov linear automorphism with $\dim E^c = 2$ is stably ergodic in the $\mathcal{C}^5$ topology.*

**Theorem 1.2.** *For $N = 4$ any pseudo-Anosov linear automorphism is stably ergodic in the $\mathcal{C}^{22}$ topology.*

As any ergodic linear automorphism of $\mathbb{T}^4$ is either Anosov or pseudo-Anosov (see Corollary A.7 of Appendix A) and any ergodic linear automorphism of $\mathbb{T}^5$ is Anosov (see Corollary A.5 of Appendix A), we get

**Corollary 1.3.** *Any ergodic linear automorphism of $\mathbb{T}^N$ is stably ergodic for $N \leq 5$.*

Thus, Corollary 1.3 solves the problem about stable ergodicity of linear automorphisms on $\mathbb{T}^N$ for $N \leq 5$. Actually, it is left the question if the differentiability assumption can be reduced. We think it may be the case, because we prove much more than ergodicity when we use the differentiability assumption. Perhaps it is as in the case of diffeomorphisms on the circle with irrational rotation number, where only a $C^2$ hypothesis leads to ergodicity with respect to Lebesgue measure. See [He] for the result about ergodicity for circle diffeomorphisms.

For $N \geq 6$ we have the same remark about the differentiability assumption. Moreover, we think that it is not too hard to weaken the hypothesis of being pseudo-Anosov to something like the linear automorphism having an Anosov part which strongly dominates a pseudo-Anosov part. Besides, using ideas like the ones in [RH] and [V], it seems that we may change the restrictive assumption $\dim E^c = 2$ by $A|_{E^c}$ being an isometry.

We point out that in [SW] Shub and Wilkinson proved that any ergodic linear automorphism can be approximated by a stably ergodic diffeomorphism, provided the dynamics on the center space is an isometry.

We finally remark that for $N$ odd, any pseudo-Anosov linear automorphisms is Anosov (see Corollary A.4 of Appendix A) so Theorem 1.1 only make sense for $N$ even. We point out that there are matrices that fit in the hypothesis of the theorem for any even $N \geq 6$. (see Proposition A.8 of Appendix A)

In order to prove ergodicity of a partially hyperbolic diffeomorphism there is a powerful property called accessibility (precise definitions are given in the next paragraph and next section). It turns out that under some hypothesis (stable under perturbations), accessibility implies ergodicity and so, stable accessibility implies stable ergodicity. Moreover, it is conjectured (see [BPSW] for instance) that accessibility is itself a stable property. However, there are systems that are stably ergodic but not partially hyperbolic. See for instance the examples of Bonatti and Viana [BV] where there the system has a dominated splitting with an expanding invariant bundle or the one studied in [T] where no invariant hyperbolic subbundle is available. Nevertheless, these examples reach a property called non-uniform hyperbolicity and moreover have some kind of accessibility property. Although the pseudo-Anosov linear automorphisms, that are not Anosov, are in fact partially hyperbolic, they are not non-uniformly hyperbolic, nor they have the accessibility property. However, they have some property called essential accessibility.

Let us concentrate on a pseudo-Anosov linear automorphisms $A$. We have that $A$ is partially hyperbolic and stably dynamically coherent (definitions will be given in the next section), so there exists a neighborhood of $A$ in the $\mathcal{C}^1$ topology where every diffeomorphism is partially hyperbolic and dynamically coherent. In the sequel we are going to shrink this neighborhood and even shrink it in the $\mathcal{C}^r$ topology. By now let us pick a diffeomorphism, say $f$, in this neighborhood. We may suppose $0$ is its fixed point. We will work mostly in the universal covering i.e. $\mathbb{R}^N$. We will call $F : \mathbb{R}^N \to \mathbb{R}^N$ the lift of $f$ that fixes the origin. We write with a tilde the points in the torus and without tilde the points in the universal covering. We denote by $W^u(x), W^s(x), W^c(x), W^{cu}(x), W^{cs}(x)$ the invariant manifolds in $\mathbb{R}^N$ i.e. the lift of the invariant manifolds to the universal covering or equivalently, the connected component of $x$ in the preimage of the invariant manifold by the covering projection. We say a path $\gamma : I \to \mathbb{T}^N$ is a $su-$path iff there exist $0 = t_0 < \cdots < t_n = 1$ such that $\gamma_i = \gamma|_{[t_i, t_{i+1}]}$ is contained in a leaf of either the $s-$foliation or the $u-$foliation. Define the equivalence relation of being $su-$accessible as always, i.e. $y$ is $su-$accessible from $x$ if there exist a $su-$path beginning at $x$ and ending at $y$. In the same way we define the relation in $\mathbb{R}^N$. Denote with $[\tilde{x}]$ the classes in the torus and with $C(x)$ the classes in the universal covering. Notice again that $C(x)$ is not a priori equal to the preimage of $[\tilde{x}]$ by the covering projection but only a connected component of it. We say that $f$ has the accessibility property if there is only one $su-$class in the torus. Moreover, we say that $f$ has the essential accessibility property if each $su-$saturated set in the torus has either null or full Lebesgue measure. By Theorem A of [PS]

**Theorem A.** *If $f \in Diff_m^2(M)$ is a center bunched, partially hyperbolic, dynamically coherent diffeomorphism with the essential accessibility property then $f$ is ergodic.*

and the properties of $A$, we only have to prove that there is a neighborhood of $A$ where each $f$ in this neighborhood has the essential accessibility property. Actually, we prove that if $C(0)$ is not trivial in certain way, then it is in fact the whole $\mathbb{R}^N$ and hence $f$ has the accessibility property. Now, if $C(0)$ is trivial and $N \geq 6$ then we use the linearizing Theorem of Arnold and Moser to get that $f$ has the essential accessibility property and is conjugated to $A$ and if $N = 4$ we adapt the result of Moser on linearization of commuting diffeomorphisms of the circle and also obtain the essential accessibility property and the conjugacy.

**Acknowledgements:** This is my thesis at IMPA under the guidance of Jacob Palis. I am very grateful to him for many valuable commentaries and all his encouragement. Also, I am indebted to Mike Shub for patiently listen the first draft of the proof and his helpful remarks. Finally, I wish to thank Enrique Pujals for several useful conversations and Raul Ures who showed me the dynamics of the pseudo-Anosov homeomorphisms of surfaces.

## 2. Preliminaries

We say that a diffeomorphism $f : M \to M$ is partially hyperbolic if there is a continuous $Df$-invariant splitting

$$TM = E_f^u \oplus E_f^c \oplus E_f^s$$

in which $E_f^s$ and $E_f^u$ are non-trivial bundles and

$$m(D^u f) > \|D^c f\| \geq m(D^c f) > \|D^s f\|$$

$$m(D^u f) > 1 > \|D^s f\|$$

where $D^\sigma f$ is the restriction of $Df$ to $E_f^\sigma$ for $\sigma = s, c$ or $u$. For a partially hyperbolic diffeomorphism, define $E_f^{cs} = E_f^c \oplus E_f^s$ and $E_f^{cu} = E_f^c \oplus E_f^u$. A partially hyperbolic diffeomorphism is dynamically coherent if the distributions $E_f^s, E_f^c, E_f^u, E_f^{cs}$ and $E_f^{cu}$ are all uniquely integrable, with the integral manifolds of $E_f^{cs}$ and $E_f^{cu}$ foliated, respectively, by the integral manifolds of $E_f^c$ and $E_f^s$ and by the integral manifolds of $E_f^c$ and $E_f^u$.

As we remark in the introduction we are going to work in $\mathbb{R}^N$. Let us recall some results. First of all the existence of the $\mathcal{F}^s, \mathcal{F}^u, \mathcal{F}^c, \mathcal{F}^{cs}, \mathcal{F}^{cu}$ foliations, the fact that each leaf is as differentiable as $f$, and depends continuously with $f$, that the center, center-stable and center-unstable foliations are plaque expansive and hence stable. Call $E^\sigma = E_A^\sigma$, $\sigma = s, u, c, cs, cu$ and $E^{su} = E^s \oplus E^u$, the invariants spaces of $A$. The same methods of construction of the invariant foliations of [HPS] let us find (see Appendix B, Proposition B.1)

$$\gamma^s : \mathbb{R}^N \times E^s \to E^{cu} \qquad \gamma^{cs} : \mathbb{R}^N \times E^{cs} \to E^u$$
$$\gamma^u : \mathbb{R}^N \times E^u \to E^{cs} \qquad \gamma^{cu} : \mathbb{R}^N \times E^{cu} \to E^s$$
$$\gamma^c : \mathbb{R}^N \times E^c \to E^{su}$$

such that calling $\gamma^\sigma(x, \cdot) = \gamma_x^\sigma, \sigma = s, u, c, cs, cu$ then

$$W^\sigma(x) = x + graph(\gamma_x^\sigma) = \{x + v + \gamma_x^\sigma(v), v \in E^\sigma\}$$

$\gamma^\sigma(x + n, v) = \gamma^\sigma(x, v)$ and $\gamma^\sigma(x, 0) = 0$. Put in $E^s, E^u, E^c$ some norm making $A|_{E^s}$ and $A^{-1}|_{E^u}$ contractions and $A|_{E^c}$ an isometry. Let us define for $v \in \mathbb{R}^N$, $|v| = |v^s| + |v^u| + |v^c|$ where $v = v^s + v^u + v^c$ with respect to $\mathbb{R}^N = E^s \oplus E^u \oplus E^c$. In the same way define for $v \in E^{cs}$, $|v| = |v^s| + |v^c|$ and the same for $E^{cu}$. It is not hard to verify that (see Appendix B)

**Lemma 2.1.** *There exist $\kappa = \kappa(f)$ such that $\kappa(f) \to 0$ as $f \xrightarrow{\mathcal{C}^1} A$ and $C > 0$ that only depends on the $\mathcal{C}^1$ size of the neighborhood of $A$ such that for $v \in E^\sigma$,*

    *(1) $|\gamma_x^\sigma(v)| \le C \log |v|$ for $\sigma = s, u, |v| \ge 2$*
    *(2) $|\gamma_x^\sigma(v)| \le \kappa$ for $\sigma = c, cs, cu$ for any $v$*
    *(3) $|(\gamma_x^u(v))^s| \le \kappa$ for any $v$*
    *(4) $|(\gamma_x^s(v))^u| \le \kappa$ for any $v$*
    *(5) $|\gamma_x^\sigma(v)| \le \kappa |v|$ for $\sigma = s, u, c, cs, cu$ for any $v$.*

We have another lemma which will be proved in Appendix B

**Lemma 2.2.** *For any $x, y \in \mathbb{R}^N$,*

    *(1) $\#W^s(x) \cap W^{cu}(y) = 1$,*
    *(2) $\#W^u(x) \cap W^{cs}(y) = 1$.*

Define

$$\pi^s : \mathbb{R}^N \to W^{cu}(0) \qquad \pi^s(x) = W^s(x) \cap W^{cu}(0),$$

$\pi^u : \mathbb{R}^N \to W^{cs}(0)$ in the same way and

$$\pi^{su} : \mathbb{R}^N \to W^c(0) \qquad \pi^{su} = \pi^s \circ \pi^u$$

Define also $j_x^\sigma : E^\sigma \to \mathbb{R}^N$, $j_x^\sigma(v) = x + v + \gamma_x^\sigma(v)$, $\sigma = s, u, c, cs, cu$ the parametrizations of the invariant manifolds.

On the other hand we have that if $f$ is $\mathcal{C}^r$ and sufficiently $\mathcal{C}^r$ near $A$ then $\mathcal{F}^s$ restricted to $W^{cs}(x)$ is a $\mathcal{C}^r$ foliation and the same holds for the $\mathcal{F}^u$ foliation (see [PSW] and Appendix B). Thus, because of Lemma 2.1, given $C > 0$ the $s$ and $u$ holonomy maps between the center manifolds of points whose center manifolds are at distance less than $C$, whenever defined, are uniformly $\mathcal{C}^r$ close to the ones of $A$. More presisly:

**Lemma 2.3.** *Given $C > 0$ and $\varepsilon > 0$ there is a neighborhood of $A$ in the $\mathcal{C}^r$ topology such that for any $f$ in this neighborhood, $x$ and $y$ with $|x - y| \leq C$, $x \in W^{cu}(y)$, calling*

$$\pi^u_{xy} : W^c(x) \to W^c(y), \qquad \pi^u_{xy}(z) = W^u(z) \cap W^c(y)$$
$$P^u_{xy} : E^c \to E^c, \qquad P^u_{xy} = (j^c_y)^{-1} \circ \pi^u_{xy} \circ j^c_x$$

*and writing*

$$P^u_{xy}(z) = z + (x - y)^c + \varphi_{xy}(z)$$

*then $\|\varphi_{xy}\|_{\mathcal{C}^r} < \varepsilon$ where we use the $\sup-norm$ in all derivatives of order less than or equal to $r$. The same holds for the $s-$holonomy.*

*Proof.* See Appendix B. $\qquad\qquad\square$

For $n \in \mathbb{Z}^N$ define

$$x_n = W^u(n) \cap W^{cs}(0) \qquad \pi^u_n : W^c(n) \to W^c(x_n)$$
$$\pi^s_n : W^c(x_n) \to W^c(0)$$

as above and

$$T_n : E^c \to E^c, \qquad T_n = (j^c_0)^{-1} \circ \pi^s_n \circ \pi^u_n \circ L_n \circ j^c_0$$

where $L_n : \mathbb{R}^N \to \mathbb{R}^N, L_n(x) = x + n$. As a consequence of the preceding lemma we have that

**Corollary 2.4.** *$T_n$ is $\mathcal{C}^r$ for all $n \in \mathbb{Z}^N$, moreover, writing*

$$T_n(z) = z + n^c + \varphi_n(z)$$

*then for any $\varepsilon > 0$ and $R > 0$ there is a neighborhood of $A$ in the $\mathcal{C}^r$ topology such that if $f$ is in this neighborhood, then $\|\varphi_n\|_{\mathcal{C}^r} < \varepsilon$ whenever $|n| \leq R$.*

Now we are going to state the linearization Theorem of Arnold and Moser. See [He].

Define for $x \in \mathbb{R}$, $\|x\| = \inf_{k \in \mathbb{Z}} |x + k|$. As usual, we say that $\alpha \in \mathbb{R}^c$ satisfy a diophantine condition with exponent $\beta$ if $\|n \cdot \alpha\| \geq \frac{c}{|n|^{c+\beta}}$ for some $c > 0$ and for any $n \in \mathbb{Z}^c$, $n \neq 0$, where $x \cdot y$ denotes the standard inner product on $\mathbb{R}^c$ and $|n| = \sum_{i=1}^c |n_i|$.

For $\alpha \in \mathbb{R}^c$, define $R_\alpha : \mathbb{R}^c \to \mathbb{R}^c$, $R_\alpha(x) = x + \alpha$. And denote with $\mathcal{C}^r_b(\mathbb{R}^c, \mathbb{R}^c)$ the set of $\mathcal{C}^r$ bounded functions.

**Theorem 2.5. KAM** *Given $\beta > 0$, $\beta \notin \mathbb{Z}$, $\alpha \in \mathbb{R}^c$ satisfying a diophantine condition with exponent $\beta$ and calling $\theta = c + \beta$, there is $V \subset \mathcal{C}^{2\theta}_b(\mathbb{R}^c, \mathbb{R}^c)$ a neighborhood of the $0$ function such that given $\varphi \in V$ satisfying $\varphi(x + n) = \varphi(x)$ for $n \in \mathbb{Z}^c$, there exist $\lambda \in \mathbb{R}^c$ and $\eta \in \mathcal{C}^\theta_b(\mathbb{R}^c, \mathbb{R}^c)$ satisfying $\eta(x + n) = \eta(x)$ for any $n \in \mathbb{Z}^c$, $\eta(0) = 0$ and such that calling $h = id + \eta$, $h$ is a diffeomorphism and calling $Q = R_\alpha + \varphi$ then $Q = R_\lambda \circ h^{-1} \circ R_\alpha \circ h$. Moreover given $\varepsilon > 0$ there is $\delta > 0$ such that if the $\mathcal{C}^{2\theta}$ size of the $\varphi$ is less than $\delta$ then the $\mathcal{C}^\theta$ size of $\eta$ and the modulus of $\lambda$ is less than $\varepsilon$.*

Let us list some properties of $A$.

**Lemma 2.6.** *For any $n \in \mathbb{Z}^N$, $n \neq 0$, and $l \in \mathbb{Z}$, $l \neq 0$, $S = \{\sum_{i=0}^{N-1} k_i A^{il} n : k_i \in \mathbb{Z} \text{ for } i = 0, \ldots N - 1\}$ is a subgroup of maximal rank.*

*Proof.* The proof follows easily from the fact that the characteristic polynomial of $A^l$ is irreducible for any nonzero $l$. See Appendix A, Lemma A.9 for more details. $\qquad\square$

Moreover, we may suppose without loss of generality that $A$ satisfy the following:

(1) $Ae_i = e_{i+1}$ for $i = 1, \ldots, N - 1$

(2) $Ae_N = -\sum_{i=0}^{N-1} p_i e_{i+1}$, $P_A(z) = \sum_{i=0}^{N} p_i z^i$ the characteristic polynomial of $A$.

Indeed, taking $n \in \mathbb{Z}^N$, $n \neq 0$, defining $L : \mathbb{R}^N \to \mathbb{R}^N$ by $L(e_i) = A^{i-1}n$ for $i = 1, \ldots N$ and taking $B = L^{-1}AL$ it is not hard to see that $B$ induces a linear automorphism and satisfies the properties listed above. Besides, given $f$ isotopic to $A$, we have its lift $F = A + \varphi$, where $\varphi$ is $\mathbb{Z}^N$−periodic and we may work with $G = B + \hat{\varphi}$ where $\hat{\varphi} = L^{-1}\varphi \circ L$ is $\mathbb{Z}^N$−periodic, and ergodicity of $G$ would imply ergodicity of $f$ as is not hard to see.

In all the paper, $C$ stands for a generic constant that only depends on the size of the neighborhood of $A$.

## 3. **Holonomies**

In this section we are going to prove some properties about the holonomies needed in the following sections.

**Proposition 3.1.** *There exists $C > 0$ that only depends on the $\mathcal{C}^1$ size of the neighborhood of $A$ and $\beta = \beta(f)$ such that $\beta(f) \to 0$ as $f \xrightarrow{\mathcal{C}^1} A$ and such that, given $x, y \in \mathbb{R}^N$, $x \in W^s(y)$, the following properties are satisfied:*

   *(1) If $d^s(x, y) \leq 2$ then $d^c(\pi^s(w), \pi^s(z)) \leq C d^c(w, z)$ for $w, z \in W^c(x)$,*
   *(2) If $d^s(x, y) \geq 1$ then $d^c(\pi^s(w), \pi^s(z)) \leq C\big(d^s(x, y)\big)^\beta d^c(w, z)$ for $w, z \in W^c(x)$,*

*And the same properties holds if $x \in W^u(y)$ interchanging $u$ and $s$.*

*Proof.* The proof of 1. is a consequence of Lemma 2.3. Let us prove 2. Take $0 < \lambda < 1$ such that $|DF|_{E^s}| < \lambda$ and $0 < \gamma = \gamma(f)$ such that $\exp(-\gamma) < |DF|_{E^c}| < \exp(\gamma)$ and we may suppose that $\gamma(f) \to 0$ as $f \xrightarrow{\mathcal{C}^1} A$. Let us take $n \geq 0$ the first integer that satisfies $d^s(F^n(y), F^n(x)) < 1$. Then we have that given $w, z \in W^c(x)$, $d^c(F^n(w), F^n(z)) \leq \exp(n\gamma)d^c(w, z)$. Now, using 1., we have that

$$d^c\big(\pi^s(F^n(w)), \pi^s(F^n(z))\big) \leq C d^c(F^n(w), F^n(z))$$

and so

$$\begin{aligned}
d^c(\pi^s(w), \pi^s(z)) &= d^c\big(F^{-n}(\pi^s(F^n(w))), F^{-n}(\pi^s(F^n(z)))\big) \\
&\leq \exp(n\gamma)d^c(\pi^s(F^n(w)), \pi^s(F^n(z))) \\
&\leq C \exp(n\gamma)d^c(F^n(w), F^n(z)) \\
&\leq C \exp(2n\gamma)d^c(w, z)
\end{aligned}$$

Let us estimate $n$. By the definition of $n$ we get that $n \leq \frac{\log d^s(y,x)}{-\log \lambda} + 1$ and so, calling $\beta = -\frac{2\gamma}{\log \lambda}$ we get

$$d^c(\pi^s(w), \pi^s(z)) \leq C \exp(2n\gamma)d^c(w, z) \leq C \exp(\gamma)\big(d^s(x, y)\big)^\beta d^c(z, w)$$

which is the desired claim.                                                                     □

**Corollary 3.2.** *There exists $C > 0$ that only depends on the neighborhood of $A$ such that for any $n \in \mathbb{Z}^N$*

   *(1) If $|n^s|, |n^u| \geq 2$ then $Lip(T_n) \leq C(|n^s||n^u|)^\beta$*
   *(2) If $|n^s| \leq 2$ and $|n^u| \geq 2$ then $Lip(T_n) \leq C|n^u|^\beta$*
   *(3) If $|n^u| \leq 2$ and $|n^s| \geq 2$ then $Lip(T_n) \leq C|n^s|^\beta$*
   *(4) If $|n^s|, |n^u| \leq 2$ then $Lip(T_n) \leq C$*

*where $\beta$ is as in Proposition 3.1.*

*Proof.* We prove the first affirmation, the others follow by the same method. We have $x_n = W^u(n) \cap W^{cs}(0)$ and $y_n = W^s(x_n) \cap W^c(0)$. So, using Proposition 3.1, we only have to estimate $d^u(n, x_n)$ and $d^s(x_n, y_n)$. Now we have that $x_n = n + v^u + \gamma_0^u(v_u) = v^{cs} + \gamma_0^{cs}(v^{cs})$ and $y_n = x_n + v^s + \gamma_{x_n}^s(v_s) = v^c + \gamma_0^c(v^c)$. So, using Lemma 2.1, we get that $|(x_n - n)^u| \leq |n^u| + \kappa$ and $|(x_n - y_n)^s| \leq |n^s| + 2\kappa$. The corollary follows from the fact $\frac{1}{C}|(x - y)^\sigma| \leq d^\sigma(x, y) \leq C|(x - y)^\sigma|$ for $\sigma = s, u, c, cs, cu$ and some constant $C > 0$ that only depends on the $\mathcal{C}^1$ size of the neighborhood of $A$. $\qquad \square$

For $L > 0$ and $x \in \mathbb{R}^N$ define $W_L^\sigma(x) = j_x^\sigma(B_L^\sigma(0))$ for $\sigma = s, u, c, cs, cu$ where $j_x^\sigma : E^\sigma \to \mathbb{R}^N$, $j_x^\sigma(v) = x + v + \gamma_x^\sigma(v)$, $\sigma = s, u, c, cs, cu$ are the parametrizations of the invariant manifolds. Moreover, call $W_L^\sigma(A) = \bigcup_{x \in A} W_L^\sigma(x)$. Given $S \subset \mathbb{Z}^N$ a subgroup of maximal rank, let us define $\mathbb{T}_S^N = \mathbb{R}^N/S$ the torus generated by the lattice $S$. Call $\nu(S) = vol(\mathbb{T}_S^N)$.

**Lemma 3.3.** *There is $b > 0$ that only depend on the size of the neighborhood of $A$ such that if we define $L(\varepsilon) = \varepsilon^{-b}$ then, given $x \in \mathbb{R}^N$ and $S \subset \mathbb{Z}^N$ subgroup of maximal rank, for $\varepsilon > 0$ small enough we have that*

$$W_\varepsilon^s(W_{L(\varepsilon)}^u(W_\varepsilon^c(x))) \cap \left( W_\varepsilon^s(W_{L(\varepsilon)}^u(W_\varepsilon^c(x))) + n \right) \neq \emptyset$$

*for some $n \in S$, $n \neq 0$.*

*Proof.* We only have to prove that there is some set $V \subset W_\varepsilon^s(W_{L(\varepsilon)}^u(W_\varepsilon^c(x)))$ such that $vol(V) > \nu(S)$. Call $W = W_\varepsilon^s(W_{L(\varepsilon)}^u(W_\varepsilon^c(x)))$. We have the following

**Claim 1.** *There is a constants $C > 0$ that only depends on the $\mathcal{C}^1$ distance of $f$ to $A$ such that for any $z \in W_{\frac{L(\varepsilon)}{2}}^u(x)$, calling $\delta = CL(\varepsilon)^{-\beta}\varepsilon$, where $\beta$ is as in Proposition 3.1, we have that $B_\delta(z) \subset W$.*

Let us left the proof of the claim until the end, and show how the lemma follows from this claim. Using the fact that $W_{\frac{L(\varepsilon)}{2}}^u(x) = j_x^u(B_{\frac{L(\varepsilon)}{2}}^u(x))$ it is not hard to see that there are points $z_1, \dots, z_n \in W_{\frac{L(\varepsilon)}{2}}^u(x)$, $n \geq C\left( L(\varepsilon)\delta^{-1} \right)^u$ where $C$ is some constant that only depends on the $\mathcal{C}^1$ size of the neighborhood of $A$ and $\dim E^u = u$ such that $W_\delta^u(z_i) \cap W_\delta^u(z_j) = \emptyset$ if $i \neq j$. Now we claim that $B_{\frac{\delta}{3}}(z_i) \cap B_{\frac{\delta}{3}}(z_j) = \emptyset$ if $i \neq j$. To prove this, write $z_i = z_j + a^u + \gamma_{z_j}^u(a^u)$ and hence we have that

$$|j_{z_j}^{-1}(z_i) - j_{z_j}^{-1}(z_j)| = |a_u| = |(z_j - z_i)^u| \leq |z_i - z_j|$$

so if the balls in $\mathbb{R}^N$ intersect, then $W_\delta^u(z_i)$ and $W_\delta^u(z_j)$ must intersect contradicting the choice of the $z_i$. Call $V = \bigcup_{i=1}^n B_{\frac{\delta}{3}}(z_i) \subset W$. Let us estimate the volume of $V$. Call $\gamma = b(u - \beta(N - u)) - (N - u)$. If $\beta$ is small enough which means if $f$ is close enough to $A$ and if we take $b$ big enough then we have that $\gamma > 0$, take for instance $b = \frac{N}{u}$, $\beta \leq \frac{u^2}{2N(N-u)}$ and so $\gamma \geq \frac{u}{2}$. So

$$vol(V) = \sum_{i=1}^n vol(B_{\frac{\delta}{3}}(z_i)) \geq Cn\delta^N \geq C\left( L(\varepsilon)\delta^{-1} \right)^u \delta^N = C\varepsilon^{-\gamma}$$

If $\varepsilon$ is small enough, as $\gamma > 0$, we get that $vol(V) > \nu(S)$. So, it is left the proof of the claim. Let us prove that for any $z \in W_{\frac{L(\varepsilon)}{2}}^u(x)$, and for any $y \in W_\varepsilon^c(x)$, $W_{\frac{2L(\varepsilon)}{3}}^u(y) \cap W^c(z) \neq \emptyset$. Call $w = W^u(y) \cap W^c(z) \neq \emptyset$ and let us show that $w \in W_{\frac{2L(\varepsilon)}{3}}^u(y)$. We have that

$$w = y + b^u + \gamma_y^u(b^u) = z + h^c + \gamma_z^c(h^c)$$
$$z = x + a^u + \gamma_x^u(a^u) \qquad y = x + r^c + \gamma_x^c(r^c)$$

and we have to estimate $|b^u|$. Now,

$$b^u = z^u - y^u + (\gamma_z^c(h^c))^u = a^u - (\gamma_x^c(r^c))^u + (\gamma_z^c(h^c))^u$$

and so $|b^u| \le |a^u| + 2\kappa \le \frac{2L(\varepsilon)}{3}$ if $\varepsilon$ is small enough which give us the intersection. Call $\pi_z^u : W^c(z) \to W^c(x)$ the unstable holonomy map. By Proposition 3.1 we have that $Lip(\pi_z^u) \le CL(\varepsilon)^\beta$ and so calling $\delta_1 = \frac{1}{C}L(\varepsilon)^{-\beta}\varepsilon$, we get $\pi_z^u(W_{\delta_1}^c(z)) \subset W_\varepsilon^c(x)$ and hence that $W_{\delta_1}^c(z) \subset W_{\frac{2L(\varepsilon)}{3}}^u(W_\varepsilon^c(x))$. Take now $y$ such that $|y - z| \le c\delta_1$ for some positive $c$ to be fixed, and define

$$w = W^s(y) \cap W^{cu}(z) \qquad r' = W^u(w) \cap W^c(x) \qquad r = W^u(w) \cap W^c(z)$$

So we want to prove that $y \in W_\varepsilon^s(w)$, $w \in W_{L(\varepsilon)}^u(r')$ and $r' \in W_\varepsilon^c(x)$. To this end, we use $r$ and so, we prove that $r \in W_{\delta_1}^c(z)$ and that $d^u(w, r)$ is small enough so that $w \in W_{L(\varepsilon)}^u(r')$.

$$y = w + a^s + \gamma_w^s(a^s) \qquad w = z + b^{cu} + \gamma_z^{cu}(b^{cu})$$

so

$$\begin{aligned} a^s &= y^s - z^s - \gamma_z^{cu}(b^{cu}) \\ b^{cu} &= (w - z)^{cu} = y^{cu} - z^{cu} - \gamma_w^s(a^s) \end{aligned}$$

and using Lemma 2.1

$$\begin{aligned} |a^s| &\le |y^s - z^s| + \kappa|b^{cu}| \\ |b^{cu}| &\le |y^{cu} - z^{cu}| + \kappa|a^s| \end{aligned}$$

which gives us

$$|a^s| + |b^{cu}| \le \frac{1}{1-\kappa}|y - z|$$

and hence

$$\begin{aligned} |a^s| &\le c_1\delta_1 \\ |b^{cu}| &\le c_1\delta_1 \end{aligned}$$

where $c_1 = \frac{c}{1-\kappa}$. Ans so, we get $y \in W_{c_1\delta_1}^s(w) \subset W_\varepsilon^s(w)$ if $\varepsilon$ is small enough. On the other hand,

$$r = z + g^c + \gamma_z^c(g^c) = w + h^u + \gamma_w^u(h^u)$$

so

$$\begin{aligned} g^c &= w^c - z^c + (\gamma_w^u(h^u))^c = (b^{cu})^c + (\gamma_w^u(h^u))^c \\ h^u &= z^u - w^u + (\gamma_z^c(g^c))^u = (b^{cu})^u + (\gamma_z^c(g^c))^u \end{aligned}$$

hence

$$\begin{aligned} |g^c| &\le c_1\delta_1 + \kappa|h^u| \\ |h^u| &\le c_1\delta_1 + \kappa|g^c| \end{aligned}$$

which gives us

$$|g^c| + |h^u| \le \frac{c_1\delta_1}{1-\kappa}$$

or

$$\begin{aligned} |g^c| &\le \frac{c}{(1-\kappa)^2}\delta_1 \\ |h^u| &\le \frac{c}{(1-\kappa)^2}\delta_1 \end{aligned}$$

So, taking $c$ sufficiently small, we get that $r \in W^c_{\delta_1}(z)$ and that $r \in W^u_{\delta_1}(w) \subset W^u_\varepsilon(w)$. Finally, as $r \in W^c_{\delta_1}(z)$, we have that $r' \in W^c_\varepsilon(x)$ and

$$r = r' + g^u + \gamma^u_{r'}(g^u)$$

we have that $|g^u| \leq \frac{2L(\varepsilon)}{3}$ and hence, as

$$w = r + t^u + \gamma^u_r(t^u) = r' + l^u + \gamma^u_{r'}(l^u)$$

we have that $l^u = t^u + g^u$ and $t^u = -h^u$, so

$$|l^u| \leq \varepsilon + \frac{2L(\varepsilon)}{3} < L(\varepsilon)$$

if $\varepsilon$ is small enough. And so $w \in W^u_{L(\varepsilon)}(r')$. $\qquad\square$

**Corollary 3.4.** *Fix $\varepsilon > 0$ as in Lemma 3.3, then $|n^u| \leq 3L(\varepsilon)$, $|n^s| \leq 4\kappa$ and $|n^c| \leq C|n^{su}|$.*

*Proof.* It follows the same spirit of the proof of Corollary 3.2. $\qquad\square$

We have another lemma.

**Lemma 3.5.** *There is $C > 0$ that only depends on the $\mathcal{C}^1$ size of the neighborhood of $A$ such that given $x \in E^c$, $n \in \mathbb{Z}^N$,*

- *(1) $|T_n(x) - (x + n^c)| \leq C \log(|n^s||n^u|) + C$, if $|n^s|, |n^u| \geq 3$*
- *(2) $|T_n(x) - (x + n^c)| \leq C \log|n^u| + C$, if $|n^u| \geq 3$ and $|n^s| \leq 3$*
- *(3) $|T_n(x) - (x + n^c)| \leq C \log|n^s| + C$, if $|n^s| \geq 3$ and $|n^u| \leq 3$*
- *(4) $|T_n(x) - (x + n^c)| \leq C$, if $|n^u|, |n^s| \leq 3$*

*Proof.* We prove the first one, the others follows in the same way. Fix $n \in \mathbb{Z}^n$, suppose $|n^s|, |n^u| \geq 3$. take $x^c \in E^c$ and $x = j^c_0(x^c)$, $y = W^u(x+n) \cap W^{cs}(0)$ and $z = W^s(y) \cap W^{cu}(0)$. Then we have that $T_n(x) = z^c$. Now

$$z - (x + n) = y - (x + n) + v^s + \gamma^s_y(v^s) = -(x + n) + w^{cu} + \gamma^{cu}_0(w^{cu})$$
$$y - (x + n) = a^u + \gamma^u_x(a^u) = -(x + n) + b^{cs} + \gamma^{cs}_0(b^{cs})$$

and hence

$$(z - (x + n))^c = (y - (x + n))^c + (\gamma^s_y(v^s))^c$$
$$(y - (x + n))^c = (\gamma^u_x(a^u))^c$$
$$a^u = -x^u - n^u + \gamma^{cs}_0(b^{cs})$$
$$v^s = -y^s + \gamma^{cu}_0(w^{cu})$$
$$y^s = x^s + n^s + (\gamma^u_x(a^u))^s$$

As $x \in W^c(0)$ we have that $|x^s|, |x^u| \leq \kappa$. So by Lemma 2.1 of Section 2 we have

$$|(z - (x + n))^c| \leq |(y - (x + n))^c| + |\gamma^s_y(v^s)|$$
$$\leq C \log|a^u| + C \log|v^s|$$
$$\leq C \log(|n^u| + 2\kappa) + C \log(\kappa + |y^s|)$$
$$\leq C \log(|n^u| + 2\kappa) + C \log(|n^s| + 3\kappa)$$

from which the result follows. $\qquad\square$

## 4. A minimal property of the system

**Theorem 4.1.** *Let $U$ be a nonempty open connected $su-$saturated subset of $\mathbb{R}^N$ and suppose there is $S \subset \mathbb{Z}^N$ a subgroup of $\mathbb{Z}^N$ of maximal rank such that $U + S = U$. Then $U = \mathbb{R}^N$.*

For the proof of the theorem we need the following:

**Proposition 4.2.** *Let $U$ be a nonempty open connected subset of $\mathbb{R}^N$ and suppose $U$ satisfies the following properties:*

*a)* $\pi_q(U) = \{0\}$ *for any $q \geq 1$,*
*b)* $U + S = U$ *for some subgroup $S \subset \mathbb{Z}^N$ of maximal rank,*

*then $U = \mathbb{R}^N$.*

*Proof.* Without loss of generality we may suppose $S = \mathbb{Z}^N$. Call $\tilde{U} = p(U)$ where $p : \mathbb{R}^N \to \mathbb{T}^N$ is the covering projection. Now, we have that $p : U \to \tilde{U}$, the restriction of $p$ to $U$, is a covering projection too. So, as $\pi_q(U) = \{0\}$ for any $q \geq 1$, we get, using Corollary 11 in Chapter 7, Section 2 of [Sp], that $\pi_q(\tilde{U}) = \{0\}$ for $q \geq 2$. Moreover, it is not hard to see that $i_\# : \pi_1(\tilde{U}) \to \pi_1(\mathbb{T}^N) = \mathbb{Z}^N$ is an isomorphism where $i_\#$ is the action of the inclusion map $i : \tilde{U} \to \mathbb{T}^N$ in the homotopy groups. Because $\tilde{U}$ is open and connected and $\pi_q(\mathbb{T}^N) = \{0\}$ for $q \geq 2$ we get that $i : \tilde{U} \to \mathbb{T}^N$ is a weak homotopy equivalence as defined after Corollary 18 in Chapter 7, Section 6 of [Sp]. As $\mathbb{T}^N$ is a CW complex, using Corollary 23 in Chapter 7, Section 6 of [Sp] we get that $i_\# : [\mathbb{T}^N; \tilde{U}] \to [\mathbb{T}^N; \mathbb{T}^N]$ is an isomorphism, where $[P; X]$ is the set of homotopy classes of maps from $P$ to $X$. Hence, there is $g : \mathbb{T}^N \to \tilde{U}$ such that $i \circ g$ is homotopic to $id : \mathbb{T}^N \to \mathbb{T}^N$. So, using degree theory, this implies that $i \circ g$ must be surjective and hence $\tilde{U} = \mathbb{T}^N$ which is equivalent to $U = \mathbb{R}^N$.   $\square$

So, we only have to prove property *a)* of the proposition. To this end, we first prove that $\pi_q(U) = \{0\}$ for $q \geq 2$ and then that $\pi_1(U) = \{0\}$. This last property is the hard one.

**Lemma 4.3.** $\pi^s : \mathbb{R}^N \to W^{cu}(0)$, $\pi^u : \mathbb{R}^N \to W^{cs}(0)$ *and* $\pi^{su} : \mathbb{R}^N \to W^c(0)$ *are fibrations (or Hurewicz fiber space) as defined at the beginning of Section 2 in Chapter 2 of* [Sp], *and so they are weak fibrations (or Serre fiber space) as defined after Corollary 4 in Chapter 7, Section 2 of* [Sp].

*Proof.* Once we prove the lemma for $\pi^s$ and $\pi^u$, the case of $\pi^{su}$ follows from Theorem 6 in Chapter 2, Section 2 of [Sp]. Let us prove then that $\pi^s$ is a fibration. Take $X$ a space, $g' : X \to \mathbb{R}^N$ and $G : X \times I \to W^{cu}(0)$ such that $G(x,0) = \pi^s \circ g'(x)$ for $x \in X$. We have to prove that there exists $G' : X \times I \to \mathbb{R}^N$ such that $G'(x,0) = g'(x)$ for $x \in X$ and $\pi^s \circ G' = G$. Define $G'(x,t) = W^s(G(x,t)) \cap W^{cu}(g'(x))$. It is not hard to see that this $G'$ makes the desired properties. The case of $\pi^u$ is completely analogous.   $\square$

**Lemma 4.4.** *Given any open and connected $s-$saturated set $E$, $\pi_q(E) = \pi_q(E \cap W^{cu}(0))$ for any $q \geq 1$. The same property holds if $E$ is $u-$saturated replacing $W^{cu}$ by $W^{cs}$. If $E$ is $su-$saturated, then $\pi_q(E) = \pi_q(E \cap W^c(0))$ for any $q \geq 1$.*

*Proof.* Since $E$ is $s-$saturated, it is not hard to see that $\pi^s|_E$ is a weak fibration and $\pi^s(E) = E \cap W^{cu}(0)$. So, take $x \in E \cap W^{cu}(0)$. As $(\pi^s)^{-1}(x) = W^s(x)$ is contractible since it is homeomorphic to $\mathbb{R}^s$ we have using Theorem 10 of Chapter 7, Section 2 of [Sp] that the following sequence

$$0 = \pi_q((\pi^s)^{-1}(x)) \overset{i_\#}{\to} \pi_q(E) \overset{\pi^s_\#}{\to} \pi_q(E \cap W^{cu}(0)) \overset{\overline{\partial}}{\to} \pi_{q-1}((\pi^s)^{-1}(x)) = 0$$

is exact and hence we get the desired result. The proof when $E$ is $u-$saturated is analogous and the case $E$ is $su-$saturated follows applying the same method to $\pi^u|_{E \cap W^{cu}(0)}$.   $\square$

**Corollary 4.5.** *Any $U$ as in Theorem 4.1 satisfies $\pi_q(U) = \{0\}$ for $q \geq 2$.*

*Proof.* By the preceding lemma we have that $\pi_q(U) = \pi_q(U \cap W^c(0))$ for any $q \geq 1$. Because $W^c(0)$ is homeomorphic to $\mathbb{R}^2$ we have that $\pi_q(U) = \pi_q(U \cap W^c(0)) = \{0\}$ for any $q \geq 2$.  □

So we want to prove that $D = U \cap W^c(0)$ is simply connected which is equivalent to prove that the complement of $D$ in the Riemann sphere is connected (looking $W^c(0)$ as $\mathbb{R}^2$), or which is equivalent, that any connected component of the complement of $D$ is not bounded.

Recall the definition of $T_n : E^c \to E^c$,

$$T_n(z) = (j_0^c)^{-1} \circ \pi_n^s \circ \pi_n^u \circ L_n \circ j_0^c$$

where $L_n : \mathbb{R}^N \to \mathbb{R}^N$, $L_n(x) = x + n$, for $n \in \mathbb{Z}^N$, $x_n = W^u(n) \cap W^{cs}(0)$ $\pi_n^u = \pi^u|_{W^c(n)}$, $\pi_n^s = \pi^s|_{W^c(x_n)}$ and $j_0^c : E^c \to \mathbb{R}^N$, $j_0^c(v) = v + \gamma_0^c(v)$ is the parametrizations of the center manifold of 0, $W^c(0)$.

Let us call $D^c = (j_0^c)^{-1}(D)$ and let us state the following proposition which solves our problem.

**Proposition 4.6.** *For any $x \in E^c$ and $\delta > 0$ there are $n \in S$, $n \neq 0$, $k \in \mathbb{Z}$, $k > 0$ and $\eta_i : [0, 1] \to E^c$, $i = 0, \ldots k - 1$ such that $\eta_i([0, 1]) \subset j_0^{-1}((C(j_0(x)) + S) \cap W^c(0))$, $\eta_i(0) \in B_\delta^c(T_{in}(x))$ and $\eta_i(1) = T_{(i+1)n}(x)$. Moreover $|T_{kn}(x) - x| \to \infty$ as $\delta \to 0$.*

Before the proof of this proposition, let us show how it solves our problem.

**Corollary 4.7.** *Any connected component of the complement of $D$ is not bounded.*

*Proof.* Take $B \subset W^c(0)$ a connected component of the complement of $D$ and call $B^c = j_0^{-1}(B)$. Take $x \in B^c$ and suppose by contradiction that $B$ is bounded. Let $R > 0$ be such that $B^c \subset B_R^c(x)$, the ball of center $x$ and radius $R$. Using the preceding proposition we have that for any $\delta > 0$ there are $n \in S$, $n \neq 0$ and $k \in \mathbb{Z}$, $k > 0$ such that $C_\delta = \bigcup_{i=0}^k \overline{B_\delta^c}(T_{in}(x)) \cup \bigcup_{i=0}^{k-1} \eta_i[0, 1]$ is connected and $|T_{kn}[0, 1] - x| \to \infty$ as $\delta \to 0$. So, for $\delta$ small enough we get that $\hat{C}_\delta$ the connected component of $C_\delta \cap \overline{B_{2R}^c}(x)$ that contains $x$ satisfies $\hat{C}_\delta \cap S_{2R}^c(x) \neq \emptyset$, where $S_{2R}^c(x)$ is the boundary of $B_{2R}^c(x)$. Then, looking at the Hausdorff space of the compact subsets of $\overline{B_{2R}^c}(x)$ we get that there is a subsequence $\delta_i \to 0$ such that $\hat{C}_{\delta_i} \to \hat{C}$ in then Hausdorff topology. Because of the properties of the Hausdorff topology, we get that $\hat{C}$ is connected, $x \in \hat{C}$, $\hat{C} \subset (E^c \setminus D^c)$, and so $\hat{C} \subset B^c$, and $\hat{C} \cap S_{2R}^c(x) \neq \emptyset$. Thus contradicting the boundedness of $B$.  □

Let us begin the proof of Proposition 4.6.

**Lemma 4.8.** *There is a constants $c > 0$ that only depends on $A$ such that calling $r = \frac{N-1}{2}$, $|n^c| \geq \frac{c}{|n|^r}$ for any $n \in \mathbb{Z}^N$, $n \neq 0$.*

*Proof.* The proof of the lemma will be carried out in Appendix A.  □

*Proof.* (**of Proposition 4.6**)

Take $\delta > 0$ and define $\varepsilon > 0$ by $\delta = \varepsilon^\gamma$, $\gamma = 1 - \beta(s + 4b)$, where $\beta$ is as in Proposition 3.1, $b$ as in Lemma 3.3, $s = rb + 1$ and $r$ as in Lemma 4.8. Moreover, we may suppose, if $f$ is sufficiently close to $A$ that $\gamma > 0$. Take $n \in S$ as in Lemma 3.3 for this $\varepsilon$. Besides, take

$\frac{\varepsilon^{-s}}{2} \leq k \leq \varepsilon^{-s}$. So, by Lemma 3.5 we have that

$$
\begin{aligned}
|T_{kn}(x) - x| & \geq |kn^c| - C \log |kn^u| - C \\
& \geq k \frac{C}{|n|^r} - bC \log 3C\varepsilon^{-1} - C \log k - C \\
& \geq \varepsilon^{-s} \frac{C}{|n^{su}|^r} - bC \log \varepsilon^{-1} - sC \log \varepsilon^{-1} - C - bC \log 3C \\
& \geq \varepsilon^{-s} C \varepsilon^{rb} - (b+s)C \log \varepsilon^{-1} - C - bC \log 3C \\
& = C\varepsilon^{-1} - (b+s)C \log \varepsilon^{-1} - C
\end{aligned}
$$

Since $\varepsilon^\gamma = \delta$, we have that $|T_{kn}(x) - x| \to \infty$ as $\delta \to 0$. Let us prove now the other part of the lemma. By Lemma 3.3, we have that

$$
W_\varepsilon^s\big(W_{L(\varepsilon)}^u(W_\varepsilon^c(j_0(x)))\big) \cap W_\varepsilon^s\big(W_{L(\varepsilon)}^u(W_\varepsilon^c(j_0(x) + n))\big) \neq \emptyset
$$

Take $z$ in this intersection. Then there are points $y, w, y', w'$ such that

$$
z \in W_\varepsilon^s(y), \ z \in W_\varepsilon^s(y'), \ y \in W_{L(\varepsilon)}^u(w), \ y' \in W_{L(\varepsilon)}^u(w' + n)
$$

and $j_0^{-1}(w) \in B_\varepsilon^c(x)$, $j_0^{-1}(w') \in B_\varepsilon^c(x)$. Now, let us define

$$
S : W^c(n) \to W^c(0), \ \ S = \pi_2^u \circ \pi^s \circ \pi_1^u,
$$

where

$$
\pi_1^u : W^c(n) \to W^c(y'), \ \ \pi^s : W^c(y') \to W^c(y), \ \ \pi_2^u : W^c(y) \to W^c(0)
$$

are the respective holonomies. By hypothesis we have that $S(w' + n) = w$. Moreover, using Proposition 3.1, we have that $Lip(S) \leq CL(\varepsilon)^{2\beta}$. Furthermore,

$$
S(j_0(x) + n) \in C(j_0(x) + n) \cap W^c(0)
$$

and hence

$$
j_0 \circ T_{in} \circ j_0^{-1} \circ S(j_0(x) + n) \in C(j_0(x) + (i+1)n) \cap W^c(0)
$$

Now, take $0 \leq i \leq k - 1$ and call $\hat{x}_{i+1} = T_{in} \circ j_0^{-1}(S(j_0(x) + n))$

$$
\begin{aligned}
d^c(\hat{x}_{i+1}, T_{in}(x)) & \leq d^c\big(\hat{x}_{i+1}, T_{in} \circ j_0^{-1} S(w' + n)\big) \\
& + d^c(T_{in}(j_0^{-1}(w)), T_{in}(x)) \\
& \leq Lip(T_{in}) Lip(j_0^{-1}) Lip(S) Lip(j_0) d^c(j_0^{-1}(w'), x) \\
& + Lip(T_{in}) d^c(j_0^{-1}(w), x) \varepsilon \\
& \leq Lip(T_{in})\big(CLip(S) + 1\big)\varepsilon
\end{aligned}
$$

Now using Corollary 3.2 we get that

$$
d^c(\hat{x}_{i+1}, T_{in}(x)) \leq (CkL(\varepsilon))^\beta \big(CL(\varepsilon)^{2\beta} + 1\big)\varepsilon \leq C\varepsilon^{1-\beta(s+3b)}
$$

If $\varepsilon$ is small enough, we get that

$$
d^c(\hat{x}_{i+1}, T_{in}(x)) < \varepsilon^{1-\beta(s+4b)} = \varepsilon^\gamma = \delta
$$

Finally, as we shall see in the next section, Lemma 5.5, there is a path $\hat{\eta}_i : [0, 1] \to W^c(0)$, $\hat{\eta}_i([0, 1]) \subset C(j_0(x) + (i+1)n)$, such that $\hat{\eta}_i(0) = j_0(\hat{x}_{i+1})$ and $\hat{\eta}_i(1) = j_0\big(T_{(i+1)n}(x)\big)$. So, taking $\eta_i = j_0^{-1} \circ \hat{\eta}_i$ we get the desired result. $\qquad \square$

As a corollary of the proof of Lemma 3.3 we have the following

**Corollary 4.9.** *Any su−saturated open subset of $\mathbb{R}^N$ has infinite volume.*

**Corollary 4.10.** *For any open $su-$saturated $U \subset \mathbb{R}^N$ and $S \subset \mathbb{Z}^N$ subgroup of maximal rank, there is $0 \neq n \in S$ such that $U \cap U + n \neq \emptyset$.*

*Proof.* $p_S : \mathbb{R}^N \to \mathbb{T}^N_S$, the covering projection to the torus generated by the lattice $S$, cannot be injective when restricted to $U$ because if it were injective we get that $vol_{\mathbb{T}^N_S}(p_S(U)) = vol_{\mathbb{R}^N}(U) = \infty$. $\qquad\square$

**Corollary 4.11.** *Any open or closed $F-$invariant $su-$saturated $U \subset \mathbb{R}^N$ satisfying $U + \mathbb{Z}^N = U$ is either empty or the whole $\mathbb{R}^N$.*

*Proof.* We prove the case $U$ is open, the case $U$ is closed follows working with the complement. Take $V \subset U$ a connected component of $U$. As $V$ is open and $su-$saturated we have by Corollary 4.10 that there is $n \in \mathbb{Z}^N$, $n \neq 0$ such that $V + n \cap V \neq \emptyset$ and so $V = V + n$ since $V + n \subset U$. Moreover, as the nonwandering set of $f$ as $\mathbb{T}^N$ we have that there are $k \in \mathbb{Z}$, $k \neq 0$ and $l \in \mathbb{Z}^N$ such that $[F^k(V) + l] \cap V \neq \emptyset$ and hence $F^k(V) + l = V$ because $F^k(V) + l \subset U$. From this, and the properties of $A$, it is not hard to see that there is a subgroup $S \subset \mathbb{Z}^N$ of maximal rank satisfying $V + S = V$. In fact, $S = \{\sum_{i=0}^{N-1} k_i A^{ik} n : k_i \in \mathbb{Z} \text{ for } i = 0, \ldots N - 1\}$. So, using Theorem 4.1 we get the corollary since $V$ is open, connected, $su-$saturated and $V + S = V$. $\qquad\square$

**Corollary 4.12.** *If $C(0)$ is open then $C(0) = \mathbb{R}^N$. And hence $f$ has the accessibility property.*

*Proof.* By Corollary 4.10 there is $n \in \mathbb{Z}^N$ such that $C(0) + n \cap C(0) \neq \emptyset$ and so $C(0) + n = C(0)$. Because $F(C(0)) = C(0)$ there is a subgroup $S \subset \mathbb{Z}^N$ of maximal rank satisfying $C(0) + S = C(0)$. Hence, as $C(0)$ is connected, using Theorem 4.1 we get that $C(0) = \mathbb{R}^N$. $\qquad\square$

## 5. Structure of the accessibility classes

In this section we are going to prove that either $C(0)$ is open, and hence the whole $\mathbb{R}^N$ by Corollary 4.12, or $\#\big(C(x) \cap W^c(0)\big) = 1$ for any $x \in \mathbb{R}^N$.

**Theorem 5.1.** *Either $C(0) = \mathbb{R}^N$ and hence $f$ has the accessibility property, or $\#\big(C(x) \cap W^c(0)\big) = 1$ for any $x \in \mathbb{R}^N$.*

The proof of the theorem essentially splits into two propositions:

**Proposition 5.2.** *For any $x \in \mathbb{R}^N$ one of the followings holds*

(1) *$C(x)$ is open*
(2) *$C(x) \cap W^c(0)$ is homeomorphic to either $S^1$, $(-1, 1)$, $[0, 1]$ or $[0, 1)$.*
(3) *$\#\big(C(x) \cap W^c(0)\big) = 1$*

*Moreover denoting with $M$ the set of points satisfying property 3. we get that $M$ is closed, $su-$saturated, $F-$invariant and $M + \mathbb{Z}^N = M$, and hence by Corollary 4.11 that $M$ is either empty or $\mathbb{R}^N$.*

**Proposition 5.3.** *In the above proposition, case 2. cannot hold for $0$ i.e. either*

(1) *$C(0)$ is open*
(2) *$\#\big(C(0) \cap W^c(0)\big) = 1$*

Before the proof of the propositions, let us prove the theorem:

*Proof.* **of Theorem 5.1** We have to prove that either $C(0) = \mathbb{R}^N$ or $M = \mathbb{R}^N$. We know that $M$ is either empty or the whole $\mathbb{R}^N$. Let us suppose that $M \neq \mathbb{R}^N$, hence $M = \emptyset$ and so, 0 must satisfy either 1. or 2. of Proposition 5.2, but by Proposition 5.3 we have that 0 must satisfy 1. and hence $C(0)$ is open and by Corollary 4.12 $C(0) = \mathbb{R}^N$. $\qquad\square$

**Lemma 5.4.** *For any $x \in \mathbb{R}^N$, $C(x) \cap W^c(0)$ is open if and only if $C(x)$ is open.*

*Proof.* If $C(x)$ is open then $C(x) \cap W^c(0)$ is open by definition of relative topology. If $C(0) \cap W^c(x)$ is open, then $(\pi^{su})^{-1}(C(x) \cap W^c(0)) = C(x)$ and hence $C(x)$ is open. $\qquad\square$

**Lemma 5.5.** *Given $x \in W^c(0)$ and $y \in C(x) \cap W^c(0)$ there is $\varepsilon_0 > 0$ and $\gamma : W^c_{\varepsilon_0}(x) \times I \to W^c(0)$ continuous such that $\gamma(x,0) = x$, $\gamma(x,1) = y$ and $\gamma(z,I) \subset C(z)$ for any $z \in W^c_{\varepsilon_0}(x)$. Where we denote by $I = [0,1]$.*

*Proof.* We first build a path in $W^c(0)$ from $x$ to $y$. Since $y \in C(x)$, there is a su path $\eta : I \to \mathbb{R}^N$ such that $\eta(0) = x$ and $\eta(1) = y$. Take $\pi^{su} \circ \eta$ and this gives the disered path. For the construction of gamma as in the lemma, just remember that the stable and unstable foliations are continuous, so if we take a point close enough to $x$, we can build a path close to $\eta$ and then project it to $W^c(0)$ as we did with $\eta$. $\qquad\square$

**Lemma 5.6.** *If $\mathrm{int}(C(x) \cap W^c(0)) \neq \emptyset$ then $C(x) \cap W^c(0)$ is open.*

*Proof.* Let $z$ and $\varepsilon > 0$ be such that $W^c_\varepsilon(z) \subset C(x) \cap W^c(0)$ and take $y \in C(x)$. Then there is $\varepsilon_0 > 0$ and $\gamma : W^c_{\varepsilon_0}(z) \times I \to W^c(0)$ continuous such that $\gamma(y,0) = y$, $\gamma(y,1) = z$ and $\gamma(w,I) \subset C(w)$ for any $w \in W^c_{\varepsilon_0}(y)$. As $\gamma(\cdot,1) = \tilde\gamma$ is continuous, $\tilde\gamma^{-1}(W^c_\varepsilon(z))$ is open and $y \in \tilde\gamma^{-1}(W^c_\varepsilon(z)) \subset C(x) \cap W^c(0)$. So $C(x) \cap W^c(0)$ is open. $\qquad\square$

By an arc we mean a homeomorphic image of $[0,1]$. In what follows let us identify $W^c(0)$ with $E^c$ for the sake of simplicity.

**Lemma 5.7.** *Suppose $C(x) \cap W^c(0)$ is not open and let $\eta_i : I \to C(x) \cap W^c(0)$ be injective $i = 1,2$. Then $\eta_1(I) \cap \eta_2(I) = \emptyset$ or $\eta_1(I) \cup \eta_2(I)$ is either an arc or a circle.*

*Proof.* Let $x$ and $\eta_i$, $i = 1,2$ be as in the lemma. Suppose that $\eta_1(I) \cap \eta_2(I) \neq \emptyset$ but that the conclusion of the lemma do not hold. We claim the following

**Claim 2.** *There are closed subintervals $I_1, I_2 \subset I$, and points $a \in I_1$, $b \in I_2$ such that $\eta_1(I_1) \cap \eta_2(I_2) = \{\eta_1(a)\} = \{\eta_2(b)\}$ and either $a \in \partial I_1$ and $b \in \mathrm{int}(I_2)$ or $a \in \mathrm{int}(I_1)$ and $b \in \partial I_2$.*

We left the proof of the claim until the end. Without loss of generality we may suppose tha $a \in \partial I_1$ and $b \in \mathrm{int}(I_2)$. Moreover, let us make a reparametrization and that sends $I_1$ to $[0,a]$ and $I_2$ to $[0,1]$. Let us keep on using the same notation $\eta_1$ and $\eta_2$ for this reparametrizations. Take $\varepsilon_0$ and $\gamma : B^c_{\varepsilon_0}(\eta_1(a)) \times I \to W^c(0)$ as in Lemma 5.5, such that $\gamma(\eta_1(a),1) = \eta_2(0)$. Given $\varepsilon_1 > 0$ small enough we can define

$$b_-(\varepsilon_1) = \sup\{s < b \text{ such that } \eta_2(s) \notin B^c_{\varepsilon_1}(\eta_1(a))\}$$

and

$$b_+(\varepsilon_1) = \inf\{s > b \text{ such that } \eta_2(s) \notin B^c_{\varepsilon_1}(\eta_1(a))\}$$

furthermore, define

$$a_-(\varepsilon_1) = \sup\{s < a \text{ such that } \eta_1(s) \notin B^c_{\varepsilon_1}(\eta_1(a))\}$$

Notice that

$$U = B^c_{\varepsilon_1}(\eta_1(a)) \smallsetminus \Big(\eta_1(a_-(\varepsilon_1),a] \cup \eta_2\big(b_-(\varepsilon_1),b_+(\varepsilon_1)\big)\Big)$$

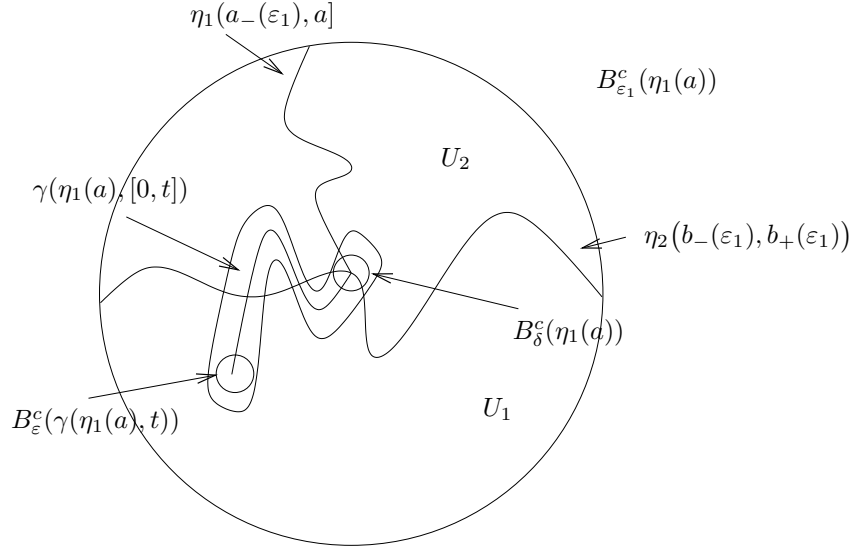

has exactly three connected components. Suppose now that there is $t > 0$ such that $\gamma(\eta_1(a), [0, t]) \subset B^c_{\varepsilon_1}(\eta_1(a))$ and $\gamma(\eta_1(a), t) \in U_1$ for $U_1$ a connected component of $U$. Take $\varepsilon > 0$ such that

$$B^c_\varepsilon(\gamma(\eta_1(a), [0, t])) \subset B^c_{\varepsilon_1}(\eta_1(a)) \qquad B^c_\varepsilon(\gamma(\eta_1(a), t)) \subset U_1$$

and take $\delta > 0$ such that if $z \in B^c_\delta(\eta_1(a))$ then $\gamma(z, s) \in B^c_\varepsilon(\gamma(\eta_1(a), s))$ for $s \in [0, t]$. Take now $U_2$ another connected component of $U$, then for any $z \in U_2 \cap B^c_\delta(\eta_1(a))$ we have that $\gamma(z, [0, t]) \subset B^c_{\varepsilon_1}(\eta_1(a))$ and $\gamma(z, t) \in U_1$. Hence, as $\gamma(z, [0, t])$ is connected we have that

$$\gamma(z, [0, t]) \cap \Big( \eta_1(a_-(\varepsilon_1), a] \cup \eta_2\big(b_-(\varepsilon_1), b_+(\varepsilon_1)\big) \Big) \neq \emptyset$$

and so $U_2 \cap B^c_\delta(\eta_1(a)) \subset C(x) \cap W^c(0)$ contradicting that it has empty interior. So we have that there is not such a $t$. This implies that defining

$$t_+ = \inf\{t > 0 \text{ such that } \gamma(\eta_1(a), t) \notin B^c_{\varepsilon_1}(\eta_1(a))\}$$

then

$$\gamma(\eta_1(a), [0, t_+)) \subset \eta_1(a_-(\varepsilon_1), a] \cup \eta_2\big(b_-(\varepsilon_1), b_+(\varepsilon_1)\big)$$

Suppose first that $\gamma(\eta_1(a), [0, t_+)) \cap \eta_2(b_-(\varepsilon_1), b] \neq \emptyset$ (the other cases will follow in a similar way). As before, we have that

$$V = B^c_{\varepsilon_1}(\eta_1(a)) \smallsetminus \big( \eta_1(a_-(\varepsilon_1), a] \cup \eta_2[b, b_+(\varepsilon_1)) \big)$$

has two connected components and that there is $0 < t' < t_+$ such that $\gamma(\eta_1(a), t') \in \eta_2(b_-(\varepsilon_1), b) \subset V_1$, where $V_1$ is one of the connected components of $V$. Because $t' < t_+$, there is $\varepsilon > 0$ such that

$$B_\varepsilon(\gamma(\eta_1(a), [0, t'])) \subset B^c_{\varepsilon_1}(\eta_1(a))$$

Now, this case follows by the same arguments as in the preceding case.

Let us prove the claim now. Call $K_i = \eta_i^{-1}(\eta_1(I) \cap \eta_2(I))$. Take $U$ a connected component of the complement of $K_1$, we have that $U \neq I$ so call $c, d$ the endpoints of $U$. Suppose $c \in K_1$ call $c' = \eta_2^{-1}(\eta_1(c))$. If $c'$ is not and endpoint of $I$, then take a point in $U$, call $I_1$ the closed interval bounded between $c$ and this point and $I_2 = I$ and we finished. If $c'$ is an endpoint of $I$, try with $d$ whenever it is in $K_1$ and work in the same way. If it do not work, then try with another connected component and with the connected components of the complement of $K_2$. If it still do not work, then it is not hard to see that we can concatenate $\eta_1$ and $\eta_2$ to build an arc or a circle. $\qquad\square$

*Proof.* **of Proposition 5.2** Suppose that 1. and 3. do not hold, then using the preceding lemma the proof follows the same spirit of the proof that the only one dimensional manifolds are the ones in 2. That $M$ is $su-$saturated, $F-$invariant and $M+\mathbb{Z}^N = M$ is almost obvious. To prove that it is closed, we prove that the complement is open. But using Lemma 5.5 it is not hard to see that the complement is in fact open. $\qquad\square$

**Corollary 5.8.** *If $f$ is sufficiently close to $A$ then $C(0) \cap W^c(0)$ is not homeomorphic to $[0,1]$.*

*Proof.* As $F(C(0) \cap W^c(0)) = C(0) \cap W^c(0)$ if $C(0) \cap W^c(0)$ were homeomorphic to $[0,1]$ then we must have that either the endpoints are fixed or permuted by $F$. As the only point fixed by $F$ is 0 and $F$ has no period two orbits we get that this is impossible. $\qquad\square$

Now, we are going to prove Proposition 5.3. So, arguing by contradiction let us suppose in the sequel that $C(0) \cap W^c(0) = \eta(J)$ where $\eta : J \to W^c(0)$, is injective, $\eta(0) = 0$ and $J = (-1,1)$ if $C(0) \cap W^c(0)$ is homeomorphic to $(-1,1)$, $J = [0,1)$ otherwise. Notice that in the case $C(0) \cap W^c(0)$ is homeomorphic to a circle, we do not have that $\eta$ is an homeomorphism. Moreover, define $H = \eta([0,1))$ and suppose, working with $f^2$ if necessary that $F(H) = H$.

We may suppose that $W^c(0)$ has the euclidean structure inherited from $E^c$. Now, we may chose $f$ close enough to $A$ in order to get the following claim.

**Claim 3.** *For $\theta \in [0,2\pi)$ define the line with slope $\theta$,*

$$l(\theta) = \{(r\cos(\theta), r\sin(\theta)) : r \geq 0\}$$

*$S(\theta)$ the sector bounded between $l(\theta)$ and $F(l(\theta))$, and $I(\theta) = \text{int}S(\theta)$. A priori there are to sectors, so, we take the one satisfying the following: there is $n \geq 2$ such that $F(I(\theta)) \cap I(\theta) = \emptyset$ and $\bigcup_{i=0}^{n} F^i(S(\theta)) = W^c(0)$. Clearly we may suppose $n$ do not depends on $\theta$, neither on $F$.*

Define $H(t) = \eta((0,t])$. As $F(H) = H$ and the only fixed point of $F$ is 0 we may suppose, working with $f^{-1}$ if necessary that $F(H(t)) \supset H(t)$.

**Lemma 5.9.** *For any $\theta$ and $t > 0$, $H(t) \cap l(\theta) \neq \emptyset$*

*Proof.* Let $\theta$ and $t > 0$ be given and suppose that $H(t) \cap l(\theta) = \emptyset$. Then we have that

$$H(t) \cap F(l(\theta)) \subset F(H(t) \cap l(\theta)) = \emptyset$$

hence $H(t) \subset I(\theta)$ or $H(t) \cap S(\theta) = \emptyset$. The first possibility can not happen because $F(I(\theta)) \cap I(\theta) = \emptyset$ and then $H(t) = H(t) \cap F(H(t)) = \emptyset$. Neither the second one because in this case we have that

$$F^{-k}(H(t)) \cap S(\theta) \subset H(t) \cap S(\theta) = \emptyset$$

hence $H(t) \cap F^k(S(\theta)) = \emptyset$ and $H(t) = H(t) \cap \bigcup_{i=0}^{n} F^i(S(\theta)) = \emptyset$. $\qquad\square$

**Corollary 5.10.** *For any $\chi : [0,\varepsilon) \to W^c(0)$ $\mathcal{C}^1$, with $\chi(0) = 0$, $\dot{\chi}(0) \neq 0$, $s > 0$ and $\delta > 0$, $\chi([0,\delta)) \cap H(s) \neq \emptyset$*

*Proof.* Take $\chi, \delta$ and $s$ as in the corollary. Call $\chi' = F \circ \chi$. As $\dot{\chi}(0) \neq \dot{\chi}'(0)$ there is $\rho > 0$ such that $\chi[0, \rho) \cap \chi'[0, \rho) = \{0\}$. Moreover, calling $CC(A, x)$ the connected component of $A$ that contains $x$, we can take $\tau$ small enough in order to get that

$$R = B_\tau^c(0) \smallsetminus \left[ CC\big(\chi[0, \rho) \cap B_\tau^c(0), 0\big) \cup CC\big(\chi'[0, \rho) \cap B_\tau^c(0), 0\big) \right]$$

has exactly two connected components. Moreover, if $\tau$ is small enough, there are $\theta_0$ and $\theta_1$ such that $l(\theta_0) \cap B_\tau^c(0) \smallsetminus \{0\}$ and $l(\theta_0) \cap B_\tau^c(0) \smallsetminus \{0\}$ do not lie in the same connected component of $R$. We may suppose that $\rho < \delta$. Take $s_0 < s$ such that $H(s_0) \subset B_\tau^c(0)$. Suppose by contradiction that $H(s_0) \cap \chi[0, \rho) = \emptyset$, then, as $F(H(s_0)) \supset H(s_0)$ we have that $H(s_0) \cap \chi'[0, \rho) = \emptyset$. By the preceding lemma, we have that $H(s_0) \cap l(\theta_0) \neq \emptyset$ and as $H(s_0)$ is connected it must lie in the same connected component in which lie $l(\theta_0) \cap B_\tau^c(0) \smallsetminus \{0\}$ thus contradicting that $H(s_0) \cap l(\theta_1) \neq \emptyset$. $\qquad \square$

Given $x \in C(0) \cap W^c(0)$ there is a $\mathcal{C}^1$ diffeomorphism $P_x : W^c(0) \to W^c(0)$ such that $P_x(0) = x$ and $P_x(z) \in C(z)$ for any $z \in W^c(0)$. To build such a diffeomorphism, take a $su-$path from $0$ to $x$ and mark the corners, then define the diffeomorphism sliding along the $s$ or $u-$foliation from the center manifold of a corner to the center manifold of the following corner. In other words, take $\gamma : [0, 1] \to \mathbb{R}^N$ call $0 = x_0, x_1, \ldots, x_n = x$ the corners of $\gamma$ enumerated by the order of $[0, 1]$, and define $\pi_0 : W^c(0) \to W^c(x_1)$ sliding along the $s-$foliation if the first leg of $\gamma$ is a $s-$path or the $u-$foliation if it is an $u-$path. Then repeat the procedure from $x_1$ to $x_2$ thus defining $\pi_1 : W^c(x_1) \to W^c(x_2)$ and so on until you reach $x_n = x$. As the holonomies we are using in the construction are at least $\mathcal{C}^1$ and the definition of accessibility class, the composition of the $\pi_i$'s give the desired $P_x$. Notice that $P_x$ is a diffeomorphism because we can make the inverse process in order to get the inverse of $P_x$.

With this $P_x$ and the corollary above we get the following

**Corollary 5.11.** *For any $t \in (0, 1)$ and any $\chi : [0, \varepsilon) \to W^c(0)$ $\mathcal{C}^1$, with $\chi(0) = \eta(t)$, $\dot{\chi}(0) \neq 0$, $s > 0$ and $\delta > 0$, $\chi((0, \delta)) \cap \eta(t - s, t + s) \neq \emptyset$*

So, let us get the contradiction to the hypothesis that $C(0) \cap W^c(0)$ is neither open nor $\{0\}$.

*Proof.* **of Proposition 5.3** Take a point $z$ in $W^c(0)$ that is not in $C(0)$. Now take the line segment from $z$ to $0$. Now take $t_0 \in (0, 1)$. Run along the line segment from $z$ to $0$ and stop the first time you touch $\eta[0, t_0]$. Suppose this point is $\eta(t_1)$. Now take the line segment from $\eta(t_1)$ to $z$ and call it $l$ then we have that $l \cap \eta[0, t_0] = \{\eta(t_1)\}$ but using the above corollary, this implies that $t_1 = t_0$ and so $\eta(t_0)$ is in the line segment from $z$ to $0$. As $t_0$ was an arbitrary point in $(0, 1)$ we have that $\eta[0, 1)$ is contained in the line segment from $z$ to $0$ thus contradicting Lemma 5.9. $\qquad \square$

## 6. **Case $C(0)$ is trivial**

We suppose in this section that $\#\big(C(x) \cap W^c(0)\big) = 1$ for any $x \in \mathbb{R}^N$.

**Lemma 6.1.** $T_n \circ T_m = T_{n+m}$ *for all $n \in \mathbb{Z}^N$.*

*Proof.* Notice that $T_n(x) = C(x + n) \cap W^c(0)$. Hence

$$
\begin{aligned}
T_n \circ T_m(x) &= C(T_m(x) + n) \cap W^c(0) \\
&= C(C(x + m) \cap W^c(0) + n) \cap W^c(0) \\
&= (C(x + m) + n) \cap W^c(0) = T_{n+m}(x)
\end{aligned}
$$

thus proving the claim. $\qquad \square$

Define the linear transformation $L : E^c \to \mathbb{R}^2$ by $L(e_1^c) = (1, 0)$ and $L(e_2^c) = (0, 1)$, call $L(n^c) = \alpha_n$ and call $P_i = L \circ T_{e_i} \circ L^{-1}$, $\tilde{Q}_n = L \circ T_n \circ L^{-1}$. Besides, take $C > 0$. We are going to chose the $\mathcal{C}^r$ neighborhood of $A$ small enough to obtain the following: There is $h : \mathbb{R}^2 \to \mathbb{R}^2$, $h = x + \eta$ such that:

(1) $h^{-1} \circ P_1 \circ h = R_{(1,0)}$,
(2) $h^{-1} \circ P_2 \circ h = R_{(0,1)}$,
(3) $h^{-1} \circ \tilde{Q}_n \circ h = Q_n$ is in some given $\mathcal{C}^r$ neighborhood of $R_{\alpha_n}$ if $|n| \leq C$, and the $\mathcal{C}^r$ neighborhood of $A$ is small enough.
(4) There is a constant $C$ such that $|\eta(z)| \leq C \log^+ |z| + C$
(5) $\eta(0) = 0$

Let us show how to build such an $h$. In the sequel, when we say that a diffeomorphism is $\mathcal{C}^r$ close to the identity, we mean that taking the $\mathcal{C}^r$ neighborhood of $A$ sufficiently small we can take the diffeomorphism as close as we want to the identity.

We have that $P_1(x, y) = (x, y) + (1, 0) + \varphi_1(x, y)$ and we can take $\|\varphi_1\|_{\mathcal{C}^r}$ as small as we want. Define $\psi : \mathbb{R} \to \mathbb{R}$, $\mathcal{C}^\infty$ such that $\psi(x) = 0$ if $x < \varepsilon$ and $\psi(x) = 1$ if $x > 1 - \varepsilon$ and we require (taking $\varepsilon$ small enough) that $\|\psi\|_{\mathcal{C}^r} \leq C$ for some fixed constant that only depends on $r$. Take $\hat{h}_1(x, y) = (x, y) + \psi(x)\varphi_1(x, y)$. If the $\mathcal{C}^r$ norm of $\varphi_1$ is sufficiently small, then $\hat{h}_1$ is a $\mathcal{C}^r$ diffeomorphism, $\mathcal{C}^r$ close to de identity. By definition we have that if $|x| \leq \varepsilon$ then

$$P_1 \circ \hat{h}_1 = \hat{h}_1 \circ R_{(1,0)} \tag{1}$$

Define now $\hat{h}(x, y) = P_1^{[x]}(\hat{h}_1(x - [x], y))$ where $[x]$ stands for the integer part of $x$. We claim that $\hat{h}|_{[-\varepsilon < x < 1 + \varepsilon]} = \hat{h}_1$, $\hat{h}$ is a $\mathcal{C}^r$ diffeomorphism and $P_1 \circ \hat{h} = \hat{h} \circ R_{(1,0)}$. The first claim is obvious if $0 \leq x < 1$. If $-\varepsilon < x < 0$ then we have that

$$\hat{h}(x, y) = P_1^{-1}(\hat{h}_1(x + 1, y)) = \hat{h}_1(x, y)$$

by (1). In the same way we get the first claim if $1 \leq x < 1 + \varepsilon$. That $\hat{h}$ is $\mathcal{C}^r$ is essentially by definition, because if we have that $(x, y)$ satisfies that $x \notin \mathbb{Z}$, then there is a neighborhood of $(x, y)$ such that $[x'] = [x]$ for any $(x', y')$ in this neighborhood and if $x \in \mathbb{Z}$ then taking the neighborhood of $(x, y)$ such that $|x - x'| < \varepsilon/2$, for $(x', y')$ in this neighborhood, and using the first part of the claim, the property follows. That $\hat{h}$ is in fact a diffeomorphism also follows in the same way, just notice that defining $n$ such that $P_1^{-n}(z, w) \in \hat{h}_1[0 \leq x < 1]$ we have that

$$\hat{h}^{-1}(z, w) = \hat{h}_1^{-1} \circ P_1^{-n}(z, w) + n$$

So we have that $\hat{h}^{-1} \circ P_1 \circ \hat{h} = R_{(1,0)}$. Define $P_2' = \hat{h}^{-1} \circ P_2 \circ \hat{h}$. By the commutativity we have that $P_2'(x + (1, 0)) = P_2'(x) + (1, 0)$ so $P_2'$ induces a diffeomorphism of the cilinder. Now taking the circle $[y = 0]$ and working us above, we can build a $\mathcal{C}^r$ diffeomorphism $h' : \mathbb{R}^2 \to \mathbb{R}^2$ with $h'(0, 0) = (0, 0)$ and such that

$$h'(x + (1, 0)) = h'(x) + (1, 0) \qquad P_2' \circ h' = h' \circ R_{(0,1)}$$

and $h'|_{[-\varepsilon < x, y < 1 + \varepsilon]}$ is $\mathcal{C}^r$ close to the identity. So taking $h = \hat{h} \circ h'$ we have that $h$ is a $\mathcal{C}^r$ diffeomorphism, $h$ restricted to some small neighborhood of the standard square is $\mathcal{C}^r$ close to the identity and

$$h^{-1} \circ P_1 \circ h = R_{(1,0)} \qquad h^{-1} \circ P_2 \circ h = R_{(0,1)}$$

Notice that the third condition must be verified only in a neighborhood of the standard square and that it is verified because of Corollary 2.4 and the fact that $h$ maybe chosen as closed to the identity as desired. That $h(0) = 0$ follows again by construction. Let us prove that $h$ satisfy condition 4. Call $x = (x_1, x_2)$ and $n = [x_1]e_1 + [x_2]e_2$, then, using the first property and Lemma 3.5, we have that

$$
\begin{aligned}
|\eta(x)| &= |h(x) - x| = |P_1^{[x_1]} \circ P_2^{[x_2]} \circ h(x - ([x_1], [x_2])) - x| \\
&\leq \left| LT_n \circ L^{-1} h(x - ([x_1], [x_2])) - x \right| \\
&\leq \|L\| \left| T_n \circ L^{-1} h(x - ([x_1], [x_2])) - \left( L^{-1} h(x - ([x_1], [x_2])) + n^c \right) \right| \\
&+ |h(x - ([x_1], [x_2])) - (x - ([x_1], [x_2]))| \\
&\leq C\|L\| \log|n| + |\eta(x - ([x_1], [x_2]))| \leq C\|L\| \log|x| + C
\end{aligned}
$$

### 6.1. Case $N \geq 6$.

**Lemma 6.2.** *If $N \geq 6$ there is $n \in \mathbb{Z}^N$ such that if we take the linear transformation $L : E^c \to \mathbb{R}^2$ defined by $L(e_1^c) = (1, 0)$, and $L(e_2^c) = (0, 1)$ and call $L(n^c) = \alpha$, then $\alpha$ satisfies a diophantine condition with exponent $\delta$ for any $\delta > 0$. Clearly, $n$ only depends on $A$.*

*Proof.* The proof of the lemma will be carried out in Appendix A. □

Now, using the KAM Theorem, we have that $Q = Q_n = R_\lambda \circ h_1^{-1} \circ R_\alpha \circ h_1$ with $\|h_1 - id\|_{\mathcal{C}^1} < \frac{1}{2}$.

**Lemma 6.3.** *Let $Q : \mathbb{R}^c \to \mathbb{R}^c$, $Q = R_\lambda \circ h_1^{-1} \circ R_\alpha \circ h_1$ and suppose $\|h_1 - id\|_{\mathcal{C}^1} < \delta$ then*

$$|Q^k(x) - (x + k\alpha)| \geq k|\lambda|(1 - \delta) + \delta(|\lambda| - 2)$$

*and*

$$|Q^k(0)| \leq k(|\lambda| + |\alpha| + 2\delta)$$

*for all $k \geq 0$*

*Proof.* Denote $h_1 = x + \varphi$ and notice that $Q(x) = x + \alpha + \lambda + \varphi(x) - \varphi(Q(x) - \lambda)$ and so

$$
\begin{aligned}
Q^k(x) &= x + k\alpha + k\lambda + \sum_{j=0}^{k-1} \varphi(Q^j(x)) - \sum_{j=0}^{k-1} \varphi(Q^{j+1}(x) - \lambda) \\
&= x + k\alpha + k\lambda + \sum_{j=1}^{k-1} [\varphi(Q^j(x)) - \varphi(Q^j(x) - \lambda)] \\
&+ \varphi(x) - \varphi(Q^k(x) - \lambda)
\end{aligned}
$$

*so*

$$\left| \sum_{j=1}^{k-1} [\varphi(Q^j(x)) - \varphi(Q^j(x) - \lambda)] + \varphi(x) - \varphi(Q^k(x) - \lambda) \right| \leq (k-1)|\lambda|\delta + 2\delta$$

and hence we get the first estimate. The second one follows easily using the same method. □

The following lemma give us another bound.

**Lemma 6.4.** *There is $C > 0$ such that for $k \geq 1$ $|\tilde{Q}_n^k(0) - k\alpha| \leq C \log k + C$ for any $k \geq 1$.*

*Proof.* It is essentially a special case of Lemma 3.5. □

Finally,

$$
\begin{aligned}
|\tilde{Q}_n^k(0) - k\alpha| &= |h(Q^k(0)) - k\alpha| = |Q^k(0) + \eta(Q^k(0)) - k\alpha| \\
&\geq k|\lambda|/2 + \frac{1}{2}(|\lambda| - 2) - |\eta(Q^k(0))| \\
&\geq k|\lambda|/2 + \frac{1}{2}(|\lambda| - 2) - C\log^+|Q^k(0)| - C \\
&\geq k|\lambda|/2 + \frac{1}{2}(|\lambda| - 2) - C\log k \\
&\quad - C\log(|\lambda| + |\alpha| + 1) - C \\
&= k|\lambda|/2 - C\log k - C
\end{aligned}
$$

and thus this implies that $\lambda = 0$.

## 6.2. Case $N = 4$.

Because Lemma 6.2 is false for $N = 4$, we need to use another argument here. What we do in this case is to show how the proof in [M] of the linearization of commuting circle diffeomorphisms applies in our case. We follow the notation of [M] in this subsection. First we need the following lemma

**Lemma 6.5.** *There exist $n_1, n_2$ such that if we take the linear transformation $L : E^c \to \mathbb{R}^2$ defined by $L(e_1^c) = (1, 0)$, and $L(e_2^c) = (0, 1)$ and call $L(n_1^c) = \alpha_1$, $L(n_2^c) = \alpha_2$, then there is a constant $c > 0$ such that*

$$
\max_{\nu = 1,2} \|\|k \cdot \alpha_\nu\|\| \geq \frac{c}{|k|^2}
$$

*for any $k \in \mathbb{Z}^2$, $k \neq 0$.*

*Proof.* The proof of the lemma will be carried out in the appendix. $\qquad\square$

So we have $\tau = 2$ in formula (1.3) of page 106 of [M]. Let us call $\phi_\nu = Q_{n_\nu}$ for $\nu = 1, 2$. Now, in page 115 of [M] formula (3.5) become

$$
V_0 = \mathcal{C}_0^\infty(\mathbb{T}^2, \mathbb{R}); \quad V_1 = \mathcal{C}_0^\infty(\mathbb{T}^2, \mathbb{R}^2); \quad V_2 = \mathcal{C}_0^\infty(\mathbb{T}^2, so(2))
$$

The operators $L, A, B, L^*, A^*, B^*$ and $M$ are defined in the same way. In Lemma 3.1 we take $\sigma = 4 + \frac{1}{30}$. Let us show how the proof of Lemma 3.1 applies in our case.

$$
\begin{aligned}
v &= \sum_{j \in \mathbb{Z}^2, \, j \neq 0} v_j e^{2\pi i j \cdot x} \\
Mv &= \sum_{j \in \mathbb{Z}^2, \, j \neq 0} \mu_j v_j e^{2\pi i j \cdot x} \\
\mu_j &= 4[\sin^2(\pi\alpha_1 \cdot j) + \sin^2(\pi\alpha_2 \cdot j)] \geq \frac{c}{|j|^2} \\
|M^{-1}v|_r &\leq C \sum_{j \in \mathbb{Z}^2, \, j \neq 0} \mu_j^{-1}|v_j||j|^r \\
&\leq C \sum_{j \in \mathbb{Z}^2, \, j \neq 0} |j|^{-(2+\frac{1}{30})}|v|_{\sigma+r} \leq C|v|_{\sigma+r}
\end{aligned}
$$

Now, in page 117, the smoothing operators are defined in the same way, changing $S^1$ by $\mathbb{T}^2$ and $\mathbb{R}$ by $\mathbb{R}^2$. The construction of the smooth solution $\tilde{u}$ is as well, defining everything componentwise. In page 118, everything works as well. The only remark is when he says that the fact $\psi_\nu$ has rotation number $\alpha_\nu$ implies $\hat{\psi}_\nu$ has a zero. In our case, we apply Lemma 3.5 which modulo changing the constants is invariant under conjugacy and thus obtain that each component of $\hat{\psi_\nu}$ has a zero, and hence in the same way we have that

$$|\hat{\psi}|_0 \leq 2|\hat{\psi} - c|_0$$

Finally in page 119, we define $\varepsilon_s = \varepsilon^{\kappa^s}$ but now, $\kappa = 1 + \frac{109}{242}$,

$$N_s = \varepsilon_s^{-\xi} \qquad l = 22 \qquad \sigma = 4 + \frac{1}{30}$$

where $\xi = \frac{2(l+5)}{l^2} = \frac{27}{242}$ and in formula (3.20) we change $|\hat{u}^{(s)}|_1 < \varepsilon_s^{1/2}$ by $|\hat{u}^{(s)}|_1 < \varepsilon_s^{4/11}$. And hence everything works and we get the conjugacy $h_1 : \mathbb{R}^2 \to \mathbb{R}^2$ whenever $|\hat{\phi}|_0$ and $|\hat{\phi}|_l$ are sufficiently small.

6.3. **End of the proof.** Call $h_2 = L^{-1} \circ h_1 \circ h^{-1} \circ L$ in either case. As the $T_n$'s form a commutative group of diffeomorphisms and the $R_{n^c}$'s acts transitively on $E^c$ we get that $h_2 \circ T_n = R_{n^c} \circ h_2$ for all $n \in \mathbb{Z}^N$. Now define $h_3 : W^c(0) \to E^c$ by $h_3 = h_2 \circ j_0^{-1}$ and $h^c : \mathbb{R}^N \to E^c$ by $h^c = h_3 \circ \pi^{su}$. We have that

$$h^c \circ L_n = h_2 \circ j_0^{-1} \circ \pi^{su} \circ L_n = h_2 \circ T_n \circ j_0^{-1} \circ \pi^{su} = R_{n^c} \circ h_2 \circ j_0^{-1} \circ \pi^{su} = R_{n^c} \circ h^c$$

moreover, $h^c|_{W^{c\sigma}(x)}$ is $C^1$ for any $x \in \mathbb{R}^N$, $\sigma = u, s$ and $y \in C(x)$ if and only if $h^c(x) = h^c(y)$. Indeed, it is not hard to see that $Lip(h^c|_{W_L^{c\sigma}(x)}) \leq C(L)$ and $Lip\big((h^c|_{W_L^c(x)})^{-1}\big) \leq C(L)$ for some constant $C(L)$ that only depends on the size of the neighborhood of $A$ and $L$. We claim that $h^c \circ F = A^c \circ h^c$. By definition, we only have to prove it in $W^c(0)$. Now,

$$h^c(F(\pi^{su}(n))) = h^c(\pi^{su}(An)) = h^c(An) = A^c n^c = A^c h^c(n) = A^c h^c(\pi^{su}(n))$$

As $\pi^{su}(\mathbb{Z}^N)$ is dense in $W^c(0)$ (this is because $h^c(\pi^{su}(\mathbb{Z}^N)) = \{n^c : n \in \mathbb{Z}^N\}$ and $h^c|_{W^c(0)}$ is a diffeomorphism), we get the desired claim. Now, denoting $F = A + \psi$ and solving the cohomological equations

$$A^s \varphi^s - \varphi^s \circ F = \psi^s \qquad \text{and} \qquad A^u \varphi^u - \varphi^u \circ F = \psi^u$$

which can be solved as in the Anosov case or as in Hartman-Grobman's Theorem, we get, as is not hard to see, that calling $h^s(x) = x^s + \varphi^s(x)$ and $h^u(x) = x^u + \varphi^u(x)$ and defining $H_1 : \mathbb{R}^N \to \mathbb{R}^N$ by $H_1 = h^s + h^u + h^c$, $H_1$ is a homeomorphism, $H_1(x + n) = H_1(x) + n$ for all $n \in \mathbb{Z}^N$ which means that $H_1$ induces a homeomorphism of the torus and that $H_1 \circ F = A \circ H_1$ and hence $f$ is conjugated to $A$. The problem now is that, as in the Anosov case, a priori $H_1$ has no regularity property other than just being continuous, hence we now define $H_2 : \mathbb{R}^N \to \mathbb{R}^N$ by $H_2(x) = x^s + x^u + h^c(x)$. $H_2$ is again a homeomorphism, and $H_2(x + n) = H_2(x) + n$ for all $n \in \mathbb{Z}^N$ and so it induces a homeomorphism of the torus. Because the properties listed above we have that $H_2(C(x)) = H_2(x) + E^{su}$. So if we prove some regularity property for $H_2$, using the fact that $x + E^{su}$, $x \in \mathbb{R}^N$ induces an ergodic foliation of the torus, we get the essential accessibility property. We claim that $H_2$ is bi-Lipschitz. To prove this claim, notice that, as $H_2$ induces a homeomorphism of the torus, it only has to be proved in a neighborhood of a fundamental domain of the torus. Moreover, we only have to prove that it is locally bi-Lipschitz by compactness. So, take $x, y$ and $\hat{x} = W^s(x) \cap W^{cu}(y)$ then $h^c(x) = h^c(\hat{x})$. Moreover, $\hat{x} = x + v^s + \gamma_x^s(v^s) = y + w^{cu} + \gamma_y^{cu}(w^{cu})$. So,

$$|\hat{x} - y| = |w^{cu} + \gamma_y^{cu}(w^{cu})| \leq (1 + \kappa)|w^{cu}|$$
$$|w^{cu}| = |(x - y)^{cu} + \gamma_x^s(v^s)| \leq |(x - y)^{cu}| + \kappa|v^s|$$
$$|v^s| = |(y - x)^s + \gamma_y^{cu}(w^{cu})| \leq |(y - x)^s| + \kappa|w^{cu}|$$
$$|w^{cu}| + |v^s| \leq \frac{1}{1 - \kappa}|x - y|$$

and $|w^{cu}| \leq \frac{1}{1-\kappa}|x-y|$. Hence

$$|\hat{x} - y| \leq \frac{1+\kappa}{1-\kappa}|x-y|$$

and so we may suppose that $x$ and $y$ are so close that $|h^c(\hat{x}) - h^c(y)| \leq C(1)|\hat{x} - y|$. As $H_2(x) = x^s + x^u + h^c(x)$ we get that $H_2$ is Lipschitz.

Let us prove now that $|H_2(x) - H_2(y)| \geq c_0|x-y|$ for some constant $c_0$. As $|H_2(x) - H_2(y)| \geq |(x-y)^{su}|$, we may suppose that $|(x-y)^{su}| \leq |(x-y)^c|$. Define $x' = W^u(\hat{x}) \cap W^c(y)$, then $x' = \hat{x} + v^u + \gamma_{\hat{x}}^u(v^u) = y + w^c + \gamma_y^c(w^c)$. So, again we may suppose $x$ and $y$ so close that

$$|h^c(x) - h^c(y)| = |h^c(x') - h^c(y)| \geq \frac{1}{C(1)}|x' - y|$$

$$|x-y| \leq 2|(x-y)^c|$$

$$|(x-y)^c| = |(w^{cu} - \gamma_x^s(v^s))^c| \leq |(\hat{x} - y)^c| + \kappa|(\hat{x} - y)^s|$$

$$|(\hat{x} - y)^c| = |w^c - (\gamma_{\hat{x}}^u(v^u))^c| \leq |(x' - y)^c| + \kappa|(x' - y)^u|$$

$$|(\hat{x} - y)^s| = |(\gamma_y^c(w^c) - \gamma_{\hat{x}}^u(v^u))^s| \leq \kappa(|(x' - y)^c| + |(x' - y)^u|)$$

Hence

$$|(x-y)^c| \leq |(x' - y)^c| + \kappa|(x' - y)^u| + \kappa^2(|(x' - y)^c| + |(x' - y)^u|) \leq 2|x' - y|$$

And so, taking $c_0 = \frac{1}{4C(1)}$ we get that $H_2^{-1}$ is Lipschitz and hence that $H_2$ is bi-Lipschitz. In fact, it can be proved that $H_2$ is a $C^1$ diffeomorphism. To do this, just notice that $h^c|_{W^c(x)}$ is a $C^1$ diffeomorphism and $h^c|_{W^\sigma(x)}$, $\sigma = s, u$ is constant for any $x \in \mathbb{R}^N$. So, it follows that the partial derivatives are continuous and hence that $h^c$ is $C^1$ and so $H_2$ is $C^1$. Working in the same way with $H_2^{-1}$ we get the desired claim.

## Appendix A. **Diophantine approximations**

In this appendix we will prove some results about diophantine approximations.

**Theorem A.1.** *Let $\alpha_i$, $i = 1, \ldots, n$ be real algebraic numbers and suppose that $1, \alpha_1, \ldots, \alpha_n$ are linearly independent over the rationals. Then, given $\delta > 0$ there is a constant $c = c(\delta, \alpha_1, \ldots, \alpha_n)$ such that for any $n+1$ integers $q_1, \ldots, q_n, p$ with $q = \max(|q_1|, \ldots, |q_n|) > 0$*

$$|q_1\alpha_1 + \cdots + q_n\alpha_n + p| \geq \frac{c}{q^{n+\delta}}$$

*Proof.* See Chapter VI, Corollary 1E of [Sc].                                    $\square$

**Proposition A.2.** *Let $P$ be a polynomial of degree $N$, with integer coefficients, irreducible over the integers. Suppose that one root of $P$ is a complex number of modulus one, say $\lambda$. Call $c_1 = 2\Re(\lambda)$ where $\Re$ stands for the real part of the number. Then, for any $Q$ with integer coefficients such that $Q(c_1) = 0$ we have that $\deg Q \geq \frac{N}{2}$.*

*Proof.* Take $Q$ such that $Q(c_1) = 0$ and $\deg Q = d$. Then, as $c_1 = \lambda + \lambda^{-1}$, we get that defining $T(x) = x^d Q(x + x^{-1})$, $\deg T = 2d$ and $T(\lambda) = 0$. As $P$ is irreducible, we get that $2d \geq N$.                                    $\square$

Let us define some tools that will be useful in what follows. For any given $\theta \in \mathbb{C}$, $|\theta| = 1$ let us denote $c_k(\theta) = 2\Re(\theta^k)$ and $a_k(\theta) = \frac{\Im(\theta^k)}{\Im(\theta)}$ where $\Im$ stands for the imaginary part of the number. We have that $a_k$ and $c_k$ satisfy the following recurrence relation:

(1) $a_0 = 0, a_1 = 1, c_0 = 2$
(2) $a_{k+1} = c_k + a_{k-1}$ for $k \geq 2$

(3) $c_k = c_1 a_k - 2a_{k-1}$ for $k \geq 1$

From this recurrence relation we get that there are polynomials with integer coefficients $R_k$ and $I_k$ that do not depend on $\theta$ such that $a_k = I_k(c_1)$ and $c_k = R_k(c_1)$, moreover, $\deg(R_k) = k$, $\deg(I_k) = k-1$ and calling $\alpha_i^k$ and $\beta_i^k$ the coefficients of $R_k$ and $I_k$ respectively, we have the following:

(1) $\alpha_k^k = 1$ and $\beta_{k-1}^k = 1$ for $k \geq 1$
(2) $\alpha_{k-2i-1}^k = 0$ and $\beta_{k-2i-2}^k = 0$ for $k \geq 1$, $0 \leq 2i \leq k-1$.

Given a polynomial $P$ with a root $\lambda$ with modulus 1 call $c_k = c_k(\lambda)$ and $a_k = a_k(\lambda)$.

**Corollary A.3.** *If $P$ is a polynomial with integer coefficients, irreducible over the integers and $\deg P$ is odd then it has no root of modulus one.*

*Proof.* Take $P$ a polynomial with integer coefficients, suppose $\deg P = 2r + 1$ and that $\lambda$ is a root of $P$ with modulus 1. Write $P(z) = \sum_{k=0}^{2r+1} p_k z^k$. Then

$$
\begin{aligned}
0 = \lambda^{-r} P(\lambda) &= \sum_{k=0}^{r} p_k \overline{\lambda}^{-k} + \sum_{k=r+1}^{2r+1} p_k \lambda^{k-r} \\
&= \sum_{k=0}^{r} p_{r-k} \overline{\lambda}^k + \sum_{k=1}^{r+1} p_{k+r} \lambda^k
\end{aligned}
$$

where $\overline{\lambda}$ is the conjugate of $\lambda$. As $\overline{\lambda}$ is also a root of $P$, we obtain, in the same way

$$
0 = \overline{\lambda}^{-r} P(\overline{\lambda}) = \sum_{k=0}^{r} p_{r-k} \lambda^k + \sum_{k=1}^{r+1} p_{k+r} \overline{\lambda}^k
$$

So, from both we obtain that

$$
\begin{aligned}
0 &= \sum_{k=0}^{r} p_{r-k}(\overline{\lambda}^k - \lambda^k) + \sum_{k=1}^{r+1} p_{k+r}(\lambda^k - \overline{\lambda}^k) \\
&= p_{2r+1}(\lambda^{r+1} - \overline{\lambda}^{r+1}) + \sum_{k=1}^{r}(p_{r+k} - p_{r-k})(\lambda^k - \overline{\lambda}^k) \\
&= 2i\left[ p_{2r+1} \mathfrak{Im}(\lambda^{r+1}) + \sum_{k=1}^{r}(p_{r+k} - p_{r-k})\mathfrak{Im}(\lambda^k) \right]
\end{aligned}
$$

So we get that

$$
0 = p_{2r+1} I_{r+1}(c_1) + \sum_{k=1}^{r}(p_{r+k} - p_{r-k}) I_k(c_1) = Q(c_1)
$$

Hence, as $\deg I_k = k - 1$ and $p_{2r+1} \neq 0$, since $\deg P = 2r + 1$, we have that $\deg Q = r$. So, using Proposition A.2 we get a contradiction, thus proving the corollary. $\square$

**Corollary A.4.** *If $N$ is odd and $A \in SL(N, \mathbb{Z})$ has irreducible characteristic polynomial then $A$ is Anosov.*

*Proof.* It is clear from the preceding corollary. $\square$

**Corollary A.5.** *Any ergodic linear automorphism of $\mathbb{T}^5$ is Anosov.*

*Proof.* Taking a power, we may suppose that $\det A = 1$. If the characteristic polynomial of $A$ where irreducible, then the result follows from the preceding corollary, so let us assume that it is reducible. Then we have that $P = LQ$ where either $\deg Q = 1, \deg L = 4$ or $\deg Q = 2, \deg L = 3$. In the first case 1 or $-1$ must be a root of $Q$ and hence of $A$ contradicting ergodicity, so we can not have this decomposition, in the second case we have that the leading coefficient of $Q$ is 1 and the independent term is $\pm 1$. Hence, if $Q$ has a root with modulus 1, it is a root of unity, contradicting ergodicity, so the roots of $Q$ do not have modulus 1. As the independent term of $L$ is also $\pm 1$, if it has a root with modulus 1, it can not be real but hence the conjugate is also a root and then, $\pm 1$ must be a root of $Q$ again contradicting ergodicity. So in this case $A$ is Anosov too.                                $\square$

**Corollary A.6.** *If $P$ is a polynomial of even degree, $\deg P = 2r$, with integers coefficients, with a root $\lambda$ of modulus one, then there is a polynomial $Q$ with integers coefficients such that $Q(c_1) = 0$ and $\deg Q = r$. Moreover, if $P$ is irreducible, denoting $P(z) = \sum_{k=0}^{2r} p_k z^k$ we have that $p_{r+k} = p_{r-k}$ for $k = 1, \ldots, r$, $Q$ is irreducible, and we can take $Q$ such that its leading coefficient equals $p_{2r}$.*

*Proof.* Here we work as in the proof of Corollary A.3.

Write $P(z) = \sum_{k=0}^{2r} p_k z^k$. Then

$$
\begin{aligned}
0 = \lambda^{-r} P(\lambda) &= p_r + \sum_{k=0}^{r-1} p_k \overline{\lambda}^{r-k} + \sum_{k=r+1}^{2r} p_k \lambda^{k-r} \\
&= p_r + \sum_{k=1}^{r} p_{r-k} \overline{\lambda}^k + \sum_{k=1}^{r} p_{k+r} \lambda^k
\end{aligned}
$$

and as $\overline{\lambda}$ is also a root of $P$, we obtain, in the same way

$$
0 = \overline{\lambda}^{-r} P(\overline{\lambda}) = p_r + \sum_{k=1}^{r} p_{r-k} \lambda^k + \sum_{k=1}^{r} p_{k+r} \overline{\lambda}^k
$$

So, from both we obtain that

$$
\begin{aligned}
0 &= 2p_r + \sum_{k=1}^{r} (p_{r+k} + p_{r-k})(\lambda^k + \overline{\lambda}^k) \\
&= 2p_r + \sum_{k=1}^{r} (p_{r+k} + p_{r-k}) 2\mathfrak{Re}(\lambda^k) \\
&= 2p_r + \sum_{k=1}^{r} (p_{r+k} + p_{r-k}) R_k(c_1) = 2Q(c_1)
\end{aligned}
$$

As $\deg R_k = k$ we get that $\deg Q \le r$ and thus we get the first part of the corollary.

For the second part we have that

$$
\begin{aligned}
0 &= \sum_{k=1}^{r} (p_{r+k} - p_{r-k})(\lambda^k - \overline{\lambda}^k) \\
&= 2i \left[ \sum_{k=1}^{r} (p_{r+k} - p_{r-k}) \mathfrak{Im}(\lambda^k) \right] \\
&= 2i \mathfrak{Im}(\lambda) \left[ \sum_{k=1}^{r} (p_{r+k} - p_{r-k}) I_k(c_1) \right] = 2i \mathfrak{Im}(\lambda) L(c_1)
\end{aligned}
$$

As $\deg I_k = k - 1$ we have that $\deg L \leq r - 1$. So, using Proposition A.2 we get that $L \equiv 0$ and so $p_{r+k} = p_{r-k}$ for $k = 1, \ldots r$. So we have that

$$Q(z) = p_r + \sum_{k=1}^{r} p_{r+k} R_k(z)$$

has degree $r$ and hence has minimal degree among the allowed, so it must be irreducible. As $\alpha_r^r = 1$, we get that the leading coefficient of $Q$ is $p_{2r}$. $\qquad \square$

**Corollary A.7.** *Any ergodic linear automorphism of $\mathbb{T}^4$ is Anosov or pseudo-Anosov.*

*Proof.* We work as in the case of $\mathbb{T}^5$. Working with $A^2$ we may suppose its determinant is 1. Suppose it is neither Anosov nor it is pseudo-Anosov. If its characteristic polynomial is irreducible then we have that $P(z) = z^4 + az^2 + 1$ but then, if $\lambda$ is a root of $P$, it must be a root of the unity, or its modulus must be different from 1. So the characteristic polynomial must be reducible, $P = LQ$, but then, $\deg L = 2$, $\deg Q = 2$ and we work as before. $\qquad \square$

**Proposition A.8.** *For any $d \geq 3$ there is a linear automorphism of $\mathbb{T}^{2d}$ in the hypothesis of Theorem 1.1.*

*Proof.* Here we will work as in Lemma A.6. For $d$ odd, define $Q(z) = z^d - 2$ and for $d$ even, define $Q(z) = z^d - d2^{d-1}z + 2$. We are going to prove that for any $d$, there are polynomials satisfying the required properties and such that calling $\lambda$ the root of modulus one and $c_1 = 2\Re(\lambda)$, then $Q(c_1) = 0$. We claim that for any $d$, $Q$ has one and only one real root $c_1$ with modulus less than or equal to 2 and it satisfies $|c_1| < 2$. For $d$ odd it is obvious as the only real root of $Q$ is $2^{1/d}$. For $d$ even, notice that $Q$ has only one minima and it is 2 and $Q(2) = 2 - (d-1)2^d < 0$ so $Q$ has exactly two real roots, one less than 2 and the other bigger than 2, moreover, $Q$ has only positive real roots and so we get the desired claim. We claim now that for any $d$, $Q$ is irreducible. Suppose by contradiction that $Q$ is reducible, $Q = LR$. We may suppose that the absolut vaule of the leading coefficient of $L$ and $R$ are both 1 and that $L(0) = 2$ and $R(0) = 1$. But it will imply (as is not hard to see) that all the coefficients of $L$ must be even, thus contradicting that its leading coefficient is $\pm 1$. Let us define $P(x) = \sum_{k=0}^{2d} p_k x^k$ where $p_{d+k} = p_{d-k}$ for $k = 1, \ldots, d$, and we want $Q(z) = p_d + \sum_{k=1}^{d} p_{d+k} R_k(z)$ where $R_k$ are defined after the proof of Proposition A.2. Now in order to find the $p_k$'s, we have to solve a $(d+1)\times(d+1)-$equation, with integer coefficients (the coefficients of the $R_k$) and as is not hard to see, it is written in a triangular form and has only ones in the diagonal. So, it has a solution in the integers, and we may choose a solution having $p_{2d} = p_0 = 1$. We claim that $P$ has only two roots with modulus one ($\lambda$ and $\overline{\lambda}$) that are not roots of the unity, it is irreducible and is not a polynomial of a power. To prove this claim, first notice that $Q$ is just the polynomial founded in Corollary A.6. So, if $P$ has another root with modulus 1 other than $\lambda$ and $\overline{\lambda}$ then $Q$ must have another real root with modulus less than or equal to 2, moreover, it must have in fact 2 roots, because if not, $\lambda = \overline{\lambda}$ and hence $\lambda = \pm 1$ and so $|c_1| = 2$. If $P$ where reducible, then there would be a polynomial $P'$ with $\deg P' < \deg P$ and $P'(\lambda) = 0$ and then we will get a polynomial $Q'$ with $\deg Q' < \deg Q$ with $Q'(c_1) = 0$ thus contradicting the irreducibility of $Q$. If it where a polynomial of a power, then it must have more that 2 roots with modulus 1, in fact, if $P(x) = T(x^n)$, $n \geq 2$, then take $\mu$ such that $\mu^n = \lambda^n$, if the only such $\mu$ are $\lambda$ and $\overline{\lambda}$ we must have that $n = 2$ and that $\lambda^2 + 1 = 0$ and thus, as $2d \geq 4$ it contradicts the irreducibility of $P$. If $\lambda$ where a root of the unity, then using the irreducibility of $P$ it is not hard to see that all the roots of $P$ must be roots of the unity and as $P$ has exactly two roots with modulus 1 then the multiplicity of $\lambda$ must be bigger that 1 thus contradicting the irreducibility. Now, defining $A$ by

(1) $Ae_i = e_{i+1}$ for $i = 1, \ldots, N-1$
(2) $Ae_N = -\sum_{i=0}^{N-1} p_i e_{i+1}$

it is not hard to see that the characteristic polynomial of $A$ is just $P$.                    □

**Lemma A.9.** *$A$ is pseudo-Anosov if and only if the characteristic polynomial of $A^l$ is irreducible for any $l \in \mathbb{Z}$, $l > 0$.*

*Proof.* If the characteristic polynomial of $A^l$ is irreducible for any $l > 0$ then the characteristic polynomial of $A$, $P_A$, is irreducible. Suppose that $P_A(x) = Q(x^n)$ for some $n \geq 2$. We have that $P_{A^n}(x^n) = P_A(x)H(x) = Q(x^n)H(x)$ where $H(x) = \det\left(\sum_{k=0}^{n-1} x^k A^{-k}\right)$. But then, it is not hard to see that $H(x) = T(x^n)$ for some polynomial $T$ and hence $P_{A^n} = QT$ thus contradicting that $P_{A^n}$ is irreducible.

Suppose now that $A$ is pseudo-Anosov but $P_{A^l}$ is reducible for some $l > 0$. Then, it is not hard to see that there is a nontrivial subgroup $S \subset \mathbb{Z}^N$ such that $A^l S = S$. Moreover there is a subgroup $R$ such that:

$$\mathbb{R}^N = [S] \oplus A[S] \oplus \cdots \oplus A^{l-1}[S] \oplus [R]$$

where $[S]$ is the subspace generated by $S$, and hence $P_A(x) = Q(x^l)T(x)$, with $Q$ and $T$ the characteristics polynomials of $A^l|_S$ and $A|_R$ respectively. As $P_A$ is irreducible and is not a polynomial of a power we get a contradiction.                    □

**Lemma A.10** (Lemma 4.8). *For any $\delta > 0$ there is a constant $c = c(A, \delta)$ such that calling $r = \frac{N}{2} - 1$ $|n^c| \geq \frac{c}{|n|^{r+\delta}}$ for any $n \in \mathbb{Z}^N$, $n \neq 0$.*

*Proof.* Call $\lambda$ the eigenvalue with modulus 1 and $e_k$ the standard basis of $\mathbb{R}^N$. Then, because of the form of $A$ we have that $A^k e_1 = e_{k+1}$ for $k = 0, \ldots, N-1$. So, given $n \in \mathbb{Z}^N$, $n = \sum_{k=0}^{N-1} n_{k+1} e_{k+1}$ we have that $n^c = \left(\sum_{k=0}^{N-1} n_{k+1}\lambda^k\right) e_1^c$. So we have that $|n^c| \geq C |\sum_{k=0}^{N-1} n_{k+1}\lambda^k|$. Now, we have using Corollary A.6 that $c_1^k = P_k(c_1)$ for any $k \geq 0$ where $P_k$ is a polynomial with integer coefficients of degree less than or equal to $\frac{N}{2} - 1$. So, we can write

$$2\Re\left(\sum_{k=0}^{N-1} n_{k+1}\lambda^k\right) = \sum_{k=0}^{N-1} n_{k+1} R_k(c_1) = \sum_{k=0}^{N/2-1} L_k^1(n) c_1^k$$

and

$$\Im\left(\sum_{k=0}^{N-1} n_{k+1}\lambda^k\right) = \Im(\lambda) \sum_{k=0}^{N-1} n_{k+1} I_k(c_1) = \Im(\lambda) \sum_{k=0}^{N/2-1} L_k^2(n) c_1^k$$

where $L_k^i$ is an homogeneous form for $k = 0, \ldots, \frac{N}{2} - 1$, $i = 1, 2$. Finally, as $c_1$ is the root of a polynomial with integers coefficients, irreducible over the integers and degree $\frac{N}{2}$ we can use Theorem A.1 and thus we get that whenever $M_1 = \max(|L_1^1(n)|, \ldots, |L_{N/2-1}^1(n)|) > 0$

$$\left|\sum_{k=0}^{N-1} n_{k+1}\lambda^k\right| \geq \left|\Re\left(\sum_{k=0}^{N-1} n_{k+1}\lambda^k\right)\right| = \frac{1}{2}\left|\sum_{k=0}^{N/2-1} L_k^1(n) c_1^k\right| \geq \frac{c}{M_1^{r+\delta}}$$

and whenever $M_2 = \max(|L_1^2(n)|, \ldots, |L_{N/2-1}^2(n)|) > 0$

$$\left|\sum_{k=0}^{N-1} n_{k+1}\lambda^k\right| \geq \left|\Im\left(\sum_{k=0}^{N-1} n_{k+1}\lambda^k\right)\right| = \frac{1}{\Im(\lambda)}\left|\sum_{k=0}^{N/2-1} L_k^2(n) c_1^k\right| \geq \frac{c}{M_2^{r+\delta}}$$

as $|L_k^i(n)| \leq C|n|$ for any $k$, $i = 1, 2$ where $C$ do not depends on $n$, the result follows whenever there is some $k \geq 1$ and $i$ such that $L_k^i(n) \neq 0$. If $L_k^i(n) = 0$ for any $k \geq 1$ and $i = 1, 2$ but $L_0^i(n) \neq 0$ the result also follows. So we must deal with the case that $L_k^i(n) = 0$ for any $k \geq 0$ and $i = 1, 2$, but this implies that $n^c = 0$ and this cannot happens, since in this case we have that $n \in E^{su}$ and the properties of $A$ implies that $E^{su} = \mathbb{R}^N$ which contradict the assumption.        □

**Lemma A.11. (Lemma 6.2)** *If $N \geq 6$ there is $n \in \mathbb{Z}^N$ such that if we take the linear transformation $L : E^c \to \mathbb{R}^2$ defined by $L(e_1^c) = (1, 0)$, and $L(e_2^c) = (0, 1)$ and call $L(n^c) = \alpha$, then $\alpha$ satisfies a diophantine condition with exponent $\delta$ for any $\delta > 0$.*

*Proof.* Take $n = e_2 + e_4$. Once we define the linear transformation that sends $e_1^c$ to $(1, 0)$, $e_2^c$ to $(0, 1)$, it is not hard to see that it sends $n^c$ to $(-c_1, c_1^2)$. Now, by Proposition A.2 we have that $1, -c_1, c_1^2$ are linearly independent over the rationals. So, using Theorem A.1 we get the lemma.        □

**Lemma A.12. (Lemma 6.5)** *If $N = 4$ there exist $n_1, n_2$ such that if we take the linear transformation $L : E^c \to \mathbb{R}^2$ defined by $L(e_1^c) = (1, 0)$, and $L(e_2^c) = (0, 1)$ and call $L(n_1^c) = \alpha_1$, $L(n_2^c) = \alpha_2$, then there is a constant $c > 0$ such that*

$$\max_{\nu = 1, 2} \|k \cdot \alpha_\nu\| \geq \frac{c}{|k|^2}$$

*for any $k \in \mathbb{Z}^2$, $k \neq 0$.*

*Proof.* Take $n_1 = (-p_1, 1 - p_2, -p_1, -1)$ and $n_2 = (1, 0, 1, 0)$ where the characteristic polynomial of $A$ is $P(z) = z^4 + p_1 z^3 + p_2 z^2 + p_1 z + 1$. Then we have that $L(n_1^c) = (c_1, 0)$ and $L(n_2^c) = (0, c_1)$. Where $c_1$ is as before for $\lambda$ the root of $P$ of modulus 1. Now we have that $Q(c_1) = c_1^2 + p_1 c_1 + p_2 - 1 = 0$ and as $Q$ is irreducible, this implies that there exists $c > 0$ such that $\|q c_1\| \geq \frac{c}{q^2}$ for any $q \in \mathbb{Z}$, $q \neq 0$. So Given $k \in \mathbb{Z}^2$, $k \neq 0$ we get that $\|k \cdot \alpha_1\| \geq \frac{c}{k_1^2} \geq \frac{c}{|k|^2}$ if $k_1 \neq 0$ and the same holds for $\alpha_2$ if $k_2 \neq 0$. As $k \neq 0$ we get the desired result.        □

## Appendix B. **Invariant Manifolds**

In this section we will show how we get the invariant foliations and how to prove regularity properties in their holonomies. We will follow [HPS] and [PSW].

**Proposition B.1.** *If $f$ is sufficiently $\mathcal{C}^r$ close to $A$ then there exists*

$$\gamma^s : \mathbb{R}^N \times E^s \to E^{cu} \qquad\qquad \gamma^{cs} : \mathbb{R}^N \times E^{cs} \to E^u$$
$$\gamma^u : \mathbb{R}^N \times E^u \to E^{cs} \qquad\qquad \gamma^{cu} : \mathbb{R}^N \times E^{cu} \to E^s$$
$$\gamma^c : \mathbb{R}^N \times E^c \to E^{su}$$

*such that calling $\gamma^\sigma(x, \cdot) = \gamma_x^\sigma$, $\sigma = s, u, c, cs, cu$ then*

$$W^\sigma(x) = x + graph(\gamma_x^\sigma) = \{x + v + \gamma_x^\sigma(v), v \in E^\sigma\}$$

*$\gamma^\sigma(x + n, v) = \gamma^\sigma(x, v)$ and $\gamma^\sigma(x, 0) = 0$. and each $\gamma^\sigma$ is continuous in the first variable and $\mathcal{C}^r$ in the second one. Moreover, $Lip(\gamma_x^\sigma) \leq \kappa$ where $\kappa = \kappa(f)$ and $\kappa(f) \to 0$ as $f \xrightarrow{\mathcal{C}^1} A$.*

*Proof.* By the invariant manifold theory, it is known that there exists $\varepsilon > 0$ and $\gamma^\sigma : \mathbb{T}^N \times E^\sigma(\varepsilon) \to E^\nu$, where $\sigma$ and $\nu$ are related in the obvious way, with all the desired regularities. So we only have to prove the existence of the global transformations, i.e. that the invariant manifolds are locally a graph is a known fact, what is new here is that they are global graphs. We are going to prove the existence of $\gamma^u$. The existence of the others follows in analogous way changing the spaces accordingly. Let us define the space

$$G = \left\{ \begin{array}{l} \gamma : E^u \to E^{cs} \text{ continuous, such that} \\ |\gamma|_* < \infty, \ Lip(\gamma) < \infty \text{ and } |\gamma|_1 < \infty \end{array} \right\}$$

where $|\gamma|_* = \sup\limits_{v \neq 0} \frac{|\gamma(v)|}{|v|}$, $Lip(\gamma)$ is the Lipschitz constant of $\gamma$ and $|\gamma|_1 = \sup\limits_{|v| \geq 2} \frac{|\gamma(v)|}{\log |v|}$. It is not hard to see that $G$ with the norm $| \cdot |_*$ is a Banach space. Define on the Banach bundle $\mathbb{T}^N \times G$, the graph transform $\Gamma$:

$$
\begin{array}{ccc}
\mathbb{T}^N \times G & \xrightarrow{\ \Gamma\ } & \mathbb{T}^N \times G \\
\downarrow & & \downarrow \\
\mathbb{T}^N & \xrightarrow{\ f\ } & \mathbb{T}^N
\end{array}
$$

defined in the following way: for $x$ and $\gamma \in G$,

$$g^u_{x,\gamma} : E^u \to E^u \qquad g^u_{x,\gamma}(w) = A^u w + \varphi^u(x + w + \gamma(w)) - \varphi^u(x)$$

and

$$\Gamma(\gamma)(x,v) = A^{cs}\gamma\big((g^u_{x,\gamma})^{-1}(v)\big) + \varphi^{cs}\Big(x + (g^u_{x,\gamma})^{-1}(v) + \gamma\big((g^u_{x,\gamma})^{-1}(v)\big)\Big) - \varphi^{cs}(x)$$

There are constants $\kappa = \kappa(f)$ satisfying $\kappa(f) \to 0$ as $f \xrightarrow{\mathcal{C}^1} A$ and $C > 0$ that only depends on the $\mathcal{C}^1$ size of the neighborhood of $A$ such that calling

$$G(\kappa, C) = \big\{ \gamma \text{ such that } |\gamma|_* \leq \kappa, \ Lip(\gamma) \leq \kappa \text{ and } |\gamma|_1 \leq C \big\}$$

$G(\kappa, C)$ is closed in $G$ and invariant under the action of the graph transform. Moreover, $\Gamma$ acts there as a contraction. So there exist a section $\eta : \mathbb{T}^N \to G(\kappa, C)$ invariant under the graph transform. Define $\gamma^u(x,v) = \eta(p(x))(v)$ where $p : \mathbb{R}^N \to \mathbb{T}^N$ is the covering projection. The continuous dependence on $f$ follows from the continuity of the invariant section in section theorems. □

**Lemma B.2. (Lemma 2.2)** *For any $x, y \in \mathbb{R}^N$,*

  *(1) $\#W^s(x) \cap W^{cu}(y) = 1$,*
  *(2) $\#W^u(x) \cap W^{cs}(y) = 1$.*

*Proof.* As always we are going to prove only the first one. Take $x, y \in \mathbb{R}^N$, to prove that they intersect we must solve the following equation:

$$x + v^s + \gamma^s_x(v^s) = y + w^{cu} + \gamma^{cu}_y(w^{cu})$$

So, call $v^s = y^s - x^s + \gamma^{cu}_y(w^{cu})$ and define $l : E^{cu} \to E^{cu}$ by

$$l(w^{cu}) = w^{cu} + x^{cu} - y^{cu} + \gamma^s_x(y^s - x^s + \gamma^{cu}_y(w^{cu})) = w^{cu} + r(w^{cu})$$

As it is not hard to see, using the preceding proposition, $Lip(r) \leq \kappa^2$, so if $\kappa < 1$ we have that $l$ is a homeomorphism and hence there exists $w^{cu}$ such that $l(w^{cu}) = 0$. It is not hard to see now that this $w^{cu}$ and $v^s = y^s - x^s + \gamma^{cu}_y(w^{cu})$ are the only ones solving the above equation. □

**Lemma B.3. (Lemma 2.1)** *There exist $\kappa = \kappa(f)$ such that $\kappa(f) \to 0$ as $f \xrightarrow{\mathcal{C}^1} A$ and $C > 0$ that only depends on the $\mathcal{C}^1$ size of the neighborhood of $A$ such that for $v \in E^\sigma$,*

  *(1) $|\gamma^\sigma_x(v)| \leq C \log |v|$ for $\sigma = s, u, |v| \geq 2$*
  *(2) $|\gamma^\sigma_x(v)| \leq \kappa$ for $\sigma = c, cs, cu$ for any $v$*
  *(3) $|(\gamma^u_x(v))^s| \leq \kappa$ for any $v$*

*(4)* $|(\gamma_x^s(v))^u| \leq \kappa$ *for any* $v$

*(5)* $|\gamma_x^\sigma(v)| \leq \kappa|v|$ *for* $\sigma = s, u, c, cs, cu$ *for any* $v$.

*Proof.* The proof of 1. and 5. follows from Proposition B.1. The proof of 3. and 4. follows from 2. as the stable and unstable manifolds subfoliates the center-stable and center-unstables manifolds. Let us prove 2. for the case of $\sigma = cu$, the other cases works as well. Denote $F = A + \psi$ and let us solve the cohomological equation

$$A^s \varphi^s - \varphi^s \circ F = \psi^s$$

Call

$$\varphi^s = \sum_{k=0}^{+\infty} (A^s)^k \psi^s \circ F^{-(k+1)}$$

We have that $\varphi^s(x + n) = \varphi^s(x)$ for any $n \in \mathbb{Z}^N$,

$$\|\varphi^s\|_0 \leq \frac{1}{1 - \|A^s\|} \|\psi^s\|_0$$

where $\| \cdot \|_0$ is the sup-norm so if $f$ is sufficiently $\mathcal{C}^0$ close to $A$ then we may suppose $\|\varphi^s\| \leq \kappa/2$. Define $h^s : \mathbb{R}^N \to E^s$ by $h^s(x) = x^s + \varphi^s(x)$. Then we have that $h^s \circ F = A^s h^s$. We claim that if $x \in W^{cu}(y)$ then $h^s(x) = h^s(y)$. Indeed

$$
\begin{aligned}
|h^s(x) - h^s(y)| &= |(A^s)^n h^s(F^{-n}(x)) - (A^s)^n h^s(F^{-n}(x))| \\
&\leq \|A^s\|^n \big( |F^{-n}(x) - F^{-n}(y)| + \kappa \big) \\
&\leq \|A^s\|^n \big( \mu^n |x - y| + \kappa \big)
\end{aligned}
$$

where $\mu = \sup \|DF^{-1}|_{E_x^{cu}}\|$ is as close to 1 as we want if $f$ is $\mathcal{C}^1$ close to $A$. And hence as we can make $\|A^s\|\mu < 1$ we have that $h^s(x) = h^s(y)$. Now, we claim that if $h^s(x) = h^s(y)$ then $W^{cu}(x) = W^{cu}(y)$. Take $x$ and $y$ such that $h^s(x) = h^s(y)$. Call $z = W^s(x) \cap W^{cu}(y)$, we claim that $z = x$. As $z \in W^{cu}(y)$, we have that

$$
\begin{aligned}
0 &= |h^s(F^n(x)) - h^s(F^n(z))| \\
&= |(F^n(x) - F^n(z))^s + (\varphi^s(F^n(x)) - \varphi^s(F^n(z)))| \\
&\geq |(F^n(x) - F^n(z))^s| - \kappa \geq C|F^n(x) - F^n(z)| - \kappa
\end{aligned}
$$

This last inequality follows because $z \in W^s(x)$. But then, letting $n \to -\infty$ we get a contradiction if $x \neq z$. So we get that $h^s(x) = h^s(y)$ if and only if $W^{cu}(x) = W^{cu}(y)$. Call now $H^s : \mathbb{R}^N \to \mathbb{R}^N$, $H^s(x) = x^{cu} + h^s(x)$. Using this last properties, it is not hard to see that $H^s(x + n) = H^s(x) + n$ and that $H^s$ is an homeomorphism, moreover we have that $H^s(W^{cu}(x)) = H^s(x) + E^{cu}$ and that $(H^s)^{-1}(y) = \hat{h}^s(y) + y^{cu}$ for some $\hat{h}^s(y) = y^s + \hat{\varphi}^s$, and so,

$$W^{cu}(x) = (H^s)^{-1}(H^s(x) + E^{cu}) = \{x + v + \varphi^s(x) + \hat{\varphi}^s(H^s(x) + v) \text{ st } v \in E^{cu}\}$$

Hence we get that $\gamma_x^{cu}(v) = \varphi^s(x) + \hat{\varphi}^s(H^s(x) + v)$, and the proof then follows from the fact that $\hat{\varphi}^s(y) = -\varphi^s((H^s)^{-1}(y))$. $\qquad\square$

**Lemma B.4. (Lemma 2.3)** *Given $C > 0$ and $\varepsilon > 0$ there is a neighborhood of $A$ in the $\mathcal{C}^r$ topology such that for any $f$ in this neighborhood, $x$ and $y$ with $|x - y| \leq C$, $x \in W^{cu}(y)$, calling*

$$\pi_{xy}^u : W^c(x) \to W^c(y), \qquad \pi_{xy}^u(z) = W^u(z) \cap W^c(y)$$

$$P_{xy}^u : E^c \to E^c, \qquad P_{xy}^u = (j_y^c)^{-1} \circ \pi_{xy}^u \circ j_x^c$$

*and writing*

$$P^u_{xy}(z) = z + (x-y)^c + \varphi_{xy}(z)$$

*then* $\|\varphi_{xy}\|_{\mathcal{C}^r} < \varepsilon$ *where we use the* $\sup-norm$ *in all derivatives of order less than or equal to* $r$. *The same holds for the* $s-holonomy$.

*Proof.* To prove this lemma, we will use the $H^s$ builded at the end of the proof of the preceding lemma and Theorem 6.7 in page 86 of [HPS]. Fix $x$ and $y$ and let us try to prove that their holonomy satisfy the required property. We omit the subindex $xy$ whenever there is no confusion. Let us try to see who is $\varphi$. Define $\phi : E^c \times E^c \to E^c$ by

$$\phi(v,w) = v + (x-y)^c - w + \left[\gamma^u\Big(j_x^c(v), \big(y + \gamma_y^c(w) - x - \gamma_x^c(v)\big)^u\Big)\right]^c$$

Then it is not hard to see that $P^u$ is the one satisfying $\phi(v, P^u(v)) = 0$. So, using the implicit function theorem, we get that if $\phi$ is $\mathcal{C}^r$ and the derivative with respect to the second variable is invertible, then $P^u$ is $\mathcal{C}^r$. Moreover, it is not hard to see that if $\phi$ is close enough in the $\mathcal{C}^r$-sup-norm to $v + (x-y) - w$, then this last property is satisfied and the $\mathcal{C}^r$-sup-norm of $\varphi$ will be small. So let us see this last property. We have to see that the $\mathcal{C}^r$-sup-norm in $v$ and $w$ of

$$(v,w) \to \left[\gamma^u\Big(x + v + \gamma_x^c(v), \big(y + \gamma_y^c(w) - x - \gamma_x^c(v)\big)^u\Big)\right]^c$$

can be taken as small as we want once it is fixed the distance between $x$ and $y$ and $f$ is close to $A$. To this end, first notice that the image of the map $(v,w) \to (y - x - \gamma_x^c(v) + \gamma_y^c(w))^u$ is contained in the ball $B^u_{(y-x)^u}(2\kappa)$. So, as it is not hard to see, if the $\mathcal{C}^r$-sup-norms of $\gamma_x^c$ and $\gamma_y^c$ are small enough and the $\mathcal{C}^r$-sup-norm of $\gamma^u|_{W^c(x) \times B^u_{(y-x)^u}(2\kappa)}$ is small enough too, then, we get the desired property. That the $\mathcal{C}^r$-sup-norms of $\gamma_x^c$ and $\gamma_y^c$ are small, follows from Lemma B.1 and Lemma B.3. Let us focus our attention on $\gamma^u$. Define $\hat{F} = H^s \circ F \circ (H^s)^{-1}$, $\hat{F} : W \to W$, where $W = \bigsqcup(p + E^{cu})$ is the disjoint union of the translate of the center-unstable space. It is a non-separable $(c + u)-$dimensional manifold. Moreover, it is not hard to see that $H^s$ looked as a diffeomorphism from $W' = \bigsqcup W^{cu}(p)$ onto $W$ is $\mathcal{C}^r$ and moreover, if $F$ is $\mathcal{C}^r$ close to $A$, then it is not hard to see that $\hat{F}$ is also $\mathcal{C}^r$ close to $A$ (looking $\hat{F}$ and $A$ as diffeomorphisms from $W$ onto $W$). Now, we build the graph transform over $W$ essentially as we did in the proof of Proposition B.1 (we change $G$ by $\hat{G}$, where $\hat{G}$ is defined as $G$ but changing $E^{cs}$ by $E^c$ and we also change $\mathbb{T}^N$ by $W$) and it turns out that it satisfies all the hypothesis of Theorem 6.7 of [HPS]. So, there is a $\mathcal{C}^r$ section $\eta_{\hat{F}} : W \to \hat{G}$ depending continuously with $\hat{F}$. Using the properties of the norm on $\hat{G}$ we get the desired property for $\gamma^u_{\hat{F}}$ and hence for $\gamma^u$ using $H^s : W' \to W$.                           $\square$

## References

[A]        Anosov, D. V. *Geodesic flows on closed Riemannian manifolds of negative curvature.* Trudy Mat. Inst. Steklov. **90** 1967.

[BPSW]     Keith Burns, Charles Pugh, Michael Shub, Amie Wilkinson, *Recent results about stable ergodicity.* Preprint.

[BV]       Christian Bonatti, Marcelo Viana, *SRB measures for partially hyperbolic systems whose central direction is mostly contracting.* Israel J. Math. **115** (2000), 157–193.

[HPS]      Morris Hirsch, Charles Pugh, Michael Shub, *Invariant manifolds.* Lecture Notes in Mathematics, Vol. 583. Springer-Verlag, Berlin-New York, 1977.

[Ha]       Paul Halmos, *On automorphisms of compact groups.* Bull. Amer. Math. Soc. **49** (1943), 619–624.

[He]       Michael Herman, *Sur la conjugaison différentiable des difféomorphismes du cercle à des rotations.* Inst. Hautes Études Sci. Publ. Math. **49** (1979), 5–233.

[M]      Jürgen Moser, *On commuting circle mappings and simultaneous Diophantine approximations.* Math. Z. **205** (1990), no. 1, 105–121.

[PS]     Charles Pugh, Michael Shub, *Stable ergodicity and julienne quasi-conformality.* J. Eur. Math. Soc. **2** (2000), no. 1, 1–52.

[PSW]    Charles Pugh, Michael Shub, Amie Wilkinson, *Hölder foliations.* Duke Math. J. **86** (1997), no. 3, 517–546.

[RH]     Jana Rodriguez Hertz, *Stable and unstable sets in a $C^0$ open and dense set of diffeomorphisms.* Publ. Mat. Urug. **8** (1999), 25–46.

[Sc]     Wolfgang Schmidt, *Diophantine approximation.* Lecture Notes in Mathematics, Vol. 785. Springer-Verlag, Berlin, 1980.

[Sp]     Edwin Spanier, *Algebraic topology.* McGraw-Hill Book Co., New York-Toronto, Ont.-London (1966)

[SW]     Michael Shub, Amie Wilkinson, *Stably ergodic approximation: two examples.* Ergodic Theory Dynam. Systems **20** (2000), no. 3, 875–893.

[T]      Ali Tahzibi, *Stably ergodic systems that are not partially hyperbolic.* Ph.D thesis at IMPA. (in preparation)

[V]      José L. Vieitez, *Expansive homeomorphisms and hyperbolic diffeomorphisms on 3-manifolds.* Erg. Th. & Dyn. Sys. **16** (1996), no.3, 591–622.

IMPA, Rio de Janeiro, Brasil
*E-mail address*: `fede@impa.br`